\documentclass[12pt]{amsart}
\usepackage{amscd}
\usepackage{latexsym}
\usepackage{amssymb}

\newtheorem{theorem}{Theorem}[section]
\newtheorem{proposition}[theorem]{Proposition}
\newtheorem{prop}[theorem]{Proposition}
\newtheorem{lemma}[theorem]{Lemma}

\newtheorem{definition}[theorem]{Definition}

\newtheorem{corollary}[theorem]{Corollary}
\newtheorem{remark}[theorem]{Remark}

\newcommand{\repthick}{{\rm Rep}_{n}(G)_{\rm thick}}
\newcommand{\chthick}{{\rm Ch}_{n}(G)_{\rm thick}}
\newcommand{\repdense}{{\rm Rep}_{n}(G)_{\rm dense}}
\newcommand{\chdense}{{\rm Ch}_{n}(G)_{\rm dense}}
\newcommand{\repmthick}{{\rm Rep}_{n}(G)_{m\mbox{-}{\rm thick}}}
\newcommand{\chmthick}{{\rm Ch}_{n}(G)_{m\mbox{-}{\rm thick}}}
\newcommand{\repmdense}{{\rm Rep}_{n}(G)_{m\mbox{-}{\rm dense}}}
\newcommand{\chmdense}{{\rm Ch}_{n}(G)_{m\mbox{-}{\rm dense}}}

\begin{document}

%%%%%%%%%%%%%%%%%%%%%%%%%%%%%%%%%%%%%%%%%%%%%%%%%%%%%%%%%%%%%%%%
%
%                       	Title 
%
%%%%%%%%%%%%%%%%%%%%%%%%%%%%%%%%%%%%%%%%%%%%%%%%%%%%%%%%%%%%%%%%

%\begin{center}
%\title{Special classes of irreducible representations I}
\title{Thick representations and dense representations I}

%\end{center}

%\bigskip

\author{Kazunori NAKAMOTO and Yasuhiro OMODA}
\address{Center for Medical Education and Sciences, Faculty of Medicine, 
University of Yamanashi}
\email{nakamoto@yamanashi.ac.jp}
\address{Department of Arts and Science, Akashi National College of Technology}
\email{omoda@akashi.ac.jp}
\thanks{The first author was partially supported by 
JSPS KAKENHI Grant Numbers JP23540044 and JP15K04814.}
%\begin{center}
%\author{ 
%Kazunori NAKAMOTO 
%\footnote{Center for Medical Education and Sciences, Faculty of Medicine, 
%University of Yamanashi}  and 
%Yasuhiro OMODA
%\footnote{Department of Arts and Science, Akashi National College of Technology} 
%}
%\end{center}

%\date{7 May, 2019 (ver. 1.4.0)}

\begin{abstract} 
We introduce special classes of irreducible representations of groups: 
thick representations and dense representations. 
Denseness implies thickness, and thickness implies 
irreducibility. We show that absolute thickness and 
absolute denseness are open conditions for representations.  
Thereby, we can construct the moduli schemes of absolutely thick representations 
and absolutely dense representations. 
We also describe several results and several examples 
on thick representations 
for developing a theory of thick representations. 
\end{abstract}



\maketitle

%%%%%%%%%%%%%%%%%%%%%%%%%%%%%%%%%%%%%%%%%%%%%%%%%%%%%%%%%%%%%%%%%
%
%  section 0.  Introduction
%
%%%%%%%%%%%%%%%%%%%%%%%%%%%%%%%%%%%%%%%%%%%%%%%%%%%%%%%%%%%%%%%%%

\section{Introduction}
We deal with special classes of irreducible representations of groups.  
First, we introduce the notion of thick representations.  
Let $G$ be a group. Let $V$ be an $n$-dimensional vector space over 
a field $k$. We say that a representation $\rho : G \to {\rm GL}(V)$ is {\it $m$-thick}  
if for any subspaces $V_1$ and $V_2$ of $V$ with $\dim V_1 = m$ and  $\dim V_2 = n - m$  
there exists $g \in G$ such that $(\rho(g) V_1) \oplus V_2 =V$. 
We also say that a representation $\rho : G \to {\rm GL}(V)$ is 
{\it thick} if $\rho$ is $m$-thick for each $0 < m < n$ 
(Definition~\ref{def:thick}).   

It may be expected that any irreducible representation is thick. 
Indeed, each irreducible representation of dimension at most $3$ is thick. 
However, it is not true for the case of dimension $n$ for $n \ge 4$. 
For example, the standard $4$-dimensional representation ${\Bbb C}^4$ of 
${\rm SO}_{4}({\Bbb C})$ is not thick (Proposition \ref{prop:so2ndensethick}).  
Hence it is a natural question 
when irreducible representations of dimension $n$ for $n \ge 4$ are thick. 

Next, we introduce another type of irreducible representations. 
We say that a representation $\rho : G \to {\rm GL}(V)$ is 
{\it $m$-dense} if the induced representation 
$(\wedge^m \rho) : G \to {\rm GL}(\wedge^m V)$ is irreducible. 
We also say that a representation $\rho : G \to {\rm GL}(V)$ is 
{\it dense} if $\rho$ is $m$-dense for each $0 < m < n$ 
(Definition~\ref{def:dense}).   
We can prove that denseness implies thickness and that thickness implies 
irreducibility (Corollary~\ref{cor:imply}). 
For example, the standard representation ${\Bbb C}^n$ of 
${\rm GL}_{n}({\Bbb C})$ is dense, and hence thick. 

The reason why we call such irreducible representations 
``thick" or ``dense" is because  the image of $\rho : G \to {\rm GL}(V)$ is 
thick or dense in ${\rm GL}(V)$, respectively. 
We imagine that if the image $\rho(G)$ gets larger in 
${\rm GL}(V)$, then $\rho$ may become thick or dense.  
Our purpose is to develop a theory of thick representations. 
Thickness is a simple, natural and essential concept in 
representation theory.  
%As an evidence that thickness is important,  
In the case of finite-dimensional representations of complex simple Lie groups, thick representations are equivalent to 
weight multiplicity free representations whose weight 
poset is a totally ordered set (\cite[Theorem~1.1]{Classification}).  
By this result, we will classify thick representations for 
complex simple Lie groups in \cite{Classification}. 
This is one of the characterization of weight multiplicity free 
representations whose weight poset 
is a totally ordered set.  

We will divide ``Thick representations and dense representations" 
into two parts: Part I and Part II, because it will be long.  
In Part I, we introduce thickness, denseness, realizable subspaces, 
and other notions on 
irreducible representations. 
We show basic results on thick representations and dense representations.  
For describing thickness, we introduce ``realizable subspaces". 
We say that a subspace $W$ of $\wedge^m V$ is {\it realizable} if 
there exist $v_1, v_2, \ldots, v_m \in V$ such that 
$0 \neq v_1 \wedge v_2 \wedge \cdots \wedge v_m \in W$ 
(Definition~\ref{def:realizable}). 
The notion of realizable subspaces is essential for   
describing criteria of thickness and the moduli of 
absolutely thick representations. 
Roughly speaking, thickness lives not in the world 
that linear algebra controls, but in the world that 
Grassmann algebra (or variety) controls. 
``Realizable subspaces" is one of keyphrases in 
Grassmann algebra. 

The main theorem of Part I is the following:  

\begin{theorem}[Theorem~\ref{th:mainth}]
Let ${\rm Rep}_n(G)$ be the representation variety of degree $n$ 
for a group $G$ over ${\mathbb Z}$.  
For $0 < m < n$, the absolutely $m$-thick representations in ${\rm Rep}_n(G)$ form an 
open subscheme of ${\rm Rep}_n(G)$. 
In particular, the absolutely thick representations in ${\rm Rep}_n(G)$ form an 
open subscheme of ${\rm Rep}_n(G)$. 
\end{theorem}
\noindent Here we say that a representation $\rho : G \to {\rm GL}(V)$ 
is {\it absolutely $m$-thick} (resp. {\it absolutely thick}) if 
$\rho\otimes_k \overline{k} : G \to {\rm GL}(V\otimes_k \overline{k})$ is 
$m$-thick (resp. thick) for an algebraic closure $\overline{k}$ of $k$.  
As a corollary of the main theorem, we can 
construct the moduli of absolutely thick representations 
(Theorems~\ref{th:modulithick}). 

In Part II, we will introduce 
$(i, j)$-thickness, $(i, j)$-denseness, and $m$-irreducibility as 
generalizations of $m$-thickness, $m$-denseness, and 
irreducibility, respectively.  
We will also describe the moduli of $4$-dimensional 
non-thick absolutely irreducible representations of 
the free group ${\rm F}_2$ of rank $2$. 
%We will deal with thickness from a viewpoint of  
%model theory, if possible. 

\smallskip 

The organization of this article is as follows: 
In \S 2, we introduce the notions of thickness and denseness. 
We describe fundamental properties of thickness and 
denseness, and a criterion for thickness. 
In \S3, we state the main theorem and 
prove the existence of the moduli schemes 
of absolutely thick representations and of absolutely dense representations. 
In \S 4, we investigate several results on  
realizable subspaces. We define the $r$-number $r(\wedge^m(n))$ and 
calculate them for small $m$ and $n$. 
In \S 5, we describe useful criteria for thickness of 
$4$-dimensional and $5$-dimensional representations. 
In \S 6, we introduce several examples of 
thick representations and dense representations for 
Lie groups. 

\bigskip 

The authors would like to express their gratitude to the referee for suggesting
several important points. The referee suggested Remarks \ref{remark:f3rep} and \ref{remark:table}, 
Proposition \ref{prop:referee}, and so on.

\section{$m$-thickness and $m$-denseness}

In this section, we introduce thickness and denseness. 
We describe fundamental properties of thickness and denseness, and a criterion 
for thickness. Proposition~\ref{prop:condofm-thick} is useful  
for verifying thickness of representations. 

\begin{definition}\label{def:thick}\rm
Let $G$ be a group. Let $V$ be an $n$-dimensional vector space over 
a field $k$. We say that a representation $\rho : G \to {\rm GL}(V)$ is {\it $m$-thick}  
if for any subspaces $V_1$ and $V_2$ of $V$ with $\dim V_1 = m$ and  $\dim V_2 = n - m$  
there exists $g \in G$ such that $(\rho(g) V_1) \oplus V_2 =V$ (or equivalently, $(\rho(g) V_1) \cap V_2 = \{ 0 \}$). 
We also say that a representation $\rho : G \to {\rm GL}(V)$ is 
{\it thick} if $\rho$ is $m$-thick for each $0 < m < n$.   
\end{definition}

\begin{remark}\rm
From the definition, any $n$-dimensional 
representations $\rho$ are always $0$-thick and $n$-thick. 
In particular, $\rho$ is thick if and only if 
$\rho$ is $m$-thick for each $0 \le m \le n$. 
\end{remark} 

\begin{definition}\label{def:dense}\rm
Let $G$ be a group. Let $V$ be an $n$-dimensional vector space over 
a field $k$. We say that a representation $\rho : G \to {\rm GL}(V)$ is 
{\it $m$-dense} if the induced representation 
$(\wedge^m \rho) : G \to {\rm GL}(\wedge^m V)$ is irreducible. 
We also say that a representation $\rho : G \to {\rm GL}(V)$ is 
{\it dense} if $\rho$ is $m$-dense for each $0 < m < n$.   
\end{definition}

\begin{remark}\rm 
For an $n$-dimensional representation $\rho : G \to {\rm GL}(V)$ over 
a field $k$, $\rho$ is always $0$-dense and $n$-dense because 
$\wedge^0 V \cong k$ and $\wedge^n V \cong k$. 
In particular, $\rho$ is dense if and only if 
$\rho$ is $m$-dense for each $0 \le m \le n$. 
\end{remark} 

\begin{lemma}\label{lemma:obvious1}
Let $\rho : G \to {\rm GL}(V)$ be an $n$-dimensional representation 
of a group $G$. 
For positive integers $i$ and $j$ with $i + j= n$, let us consider the $G$-equivariant perfect 
pairing $\wedge^i V \otimes \wedge^j V \stackrel{\wedge}{\longrightarrow} \wedge^{n} V \cong k$. 
For a $G$-invariant subspace $W$ of $\wedge^i V$, put $W^{\perp} := \{ y \in \wedge^j V 
\mid  x\wedge y = 0 \mbox{ for any } x \in W \}$. Then $W^{\perp}$ is a $G$-invariant subspace 
of $\wedge^j V$. 
In particular, $\wedge^i V$ is irreducible if and only if so is $\wedge^j V$.  
\end{lemma}

{\it Proof}. 
For $y \in W^{\perp}$, we have $x \wedge gy = g(g^{-1}x \wedge y) = 0$ for 
$x \in W$ and $g \in G$. Hence $W^{\perp}$ is $G$-invariant. 
The correspondence $W \mapsto W^{\perp}$ gives a bijection 
between the $G$-invariant subspaces of $\wedge^i V$ and $\wedge^j V$. 
Therefore $\wedge^i V$ is irreducible if and only if 
so is $\wedge^j V$.  
\qed 

\begin{proposition}\label{prop:dual} 
Let $\rho : G \to {\rm GL}(V)$ be an $n$-dimensional representation of a group $G$. 
For each $0 < m < n$, $\rho$ is $m$-thick (resp. $m$-dense) if and only if 
$\rho$ is $(n-m)$-thick (resp. $(n-m)$-dense). 
\end{proposition}

{\it Proof}. 
It is obvious that $m$-thickness and  $(n-m)$-thickness are equivalent. 
By using Lemma \ref{lemma:obvious1}, we see that $m$-denseness 
and $(n-m)$-denseness are 
equivalent.  
\qed

\begin{proposition}\label{prop:imp}
For any $n$-dimensional representations $\rho : G \to {\rm GL}(V)$, the following implications hold 
for $0 < m < n$:  
\begin{eqnarray} 
\begin{array}{ccccc} 
\mbox{$m$-dense} & \Longrightarrow & \mbox{$m$-thick}  & &  \\ 
 & & \Downarrow & & \\ 
\mbox{$1$-dense} & \Longleftrightarrow & \mbox{$1$-thick} & \Longleftrightarrow & \mbox{irreducible}. 
\end{array}  
\end{eqnarray} 
\end{proposition}

{\it Proof}. 
It suffices to prove ``$m\mbox{-dense} \Rightarrow m\mbox{-thick}$'', ``$\mbox{irreducible} \Rightarrow 
1\mbox{-dense}$'', and ``$m\mbox{-thick} \Rightarrow \mbox{irreducible}$''.   
First, we show ``$m\mbox{-dense} \Rightarrow m\mbox{-thick}$''. 
Assume that $\rho : G \to {\rm GL}(V)$ is $m$-dense. Let $V_1$ and $V_2$ be 
vector subspaces of $V$ with $\dim V_1 =m$ and $\dim V_2 = n-m$. 
The canonical homomorphism 
$\wedge^m V \otimes \wedge^{n-m} V \to \wedge^{n} V \cong k$ is 
a perfect pairing and $G$-equivariant. 
Let us take a basis $\{ e_1, e_2, \ldots, e_m \}$ of $V_1$ and a basis 
of $\{ f_1, f_2, \ldots, f_{n-m} \}$ of $V_2$. Because of irreducibility of 
$\wedge^m V$, the vectors $\{ (\wedge^m \rho) (g) (e_1 \wedge e_2 \wedge \cdots \wedge e_m) 
\mid g \in G \}$ 
span the vector space $\wedge^m V$.   
Hence there exists $g \in G$ such that 
$(\wedge^m \rho)(g) (e_1 \wedge e_2 \wedge \cdots \wedge e_m) \wedge 
(f_1 \wedge f_2 \wedge \cdots \wedge f_{n-m}) \neq 0$. 
This implies that $(\rho(g) V_1) \oplus V_2 = V$. Therefore $\rho$ is $m$-thick. 
 
%Next, we show that $1$-denseness, $1$-thickness, and irreducibility are equivalent. 
%It follows from the discussion above 
%that $1$-denseness implies $1$-thickness.    
Next, we show ``$\mbox{irreducible} \Rightarrow 
1\mbox{-dense}$''. This follows easily from the definition. 
%So we show that $1$-thickness implies irreducibility. 
%If $\rho$ is not irreducible, then there exists a non-trivial 
%$G$-invariant subspace $V'$ of $V$. 
%Let us take a $1$-dimensional subspace $V_1$ of $V'$ and 
%an $(n-1)$-dimensional subspace $V_2$ of $V$ such that $V' \subseteq V_2$.  
%Then for any $g \in G$ the intersection $(\rho(g) V_1) \cap V_2 \supseteq (\rho(g) V_1) \cap V' 
%= \rho(g) V_1 \neq 0$. Hence $(\rho(g) V_1) + V_2 \neq V$. 
%Therefore $\rho$ is not $1$-thick, which shows that $1$-thickness implies irreducibility. 
Finally, we show ``$m\mbox{-thick} \Rightarrow \mbox{irreducible}$''. 
Assume that $\rho$ is not irreducible. There exists a non-trivial $G$-invariant subspace 
$V'$ of $V$. Set $\ell := \dim V'$. 
Then we only need to choose suitable subspaces $V_1, V_2$ of $V$ such that $\dim V_1=m$, 
$\dim V_2 = n-m$ and $(\rho(g) V_1) + V_2 \neq V$ for any $g \in G$. 
This implies $\rho$ is not $m$-thick, which completes the proof. 
For the proof, we consider the following three cases: $\ell \le \min(m, n-m)$, 
$\ell \ge \max(m, n-m)$, and $\min(m, n-m) < \ell < \max(m, n-m)$. 
If $\ell \le \min(m, n-m)$, then 
let us take subspaces $V_1, V_2$ of $V$ such that 
$V' \subseteq V_1$ and $V' \subseteq V_2$. 
Since $\rho(g) V_1 \supseteq \rho(g) V' = V'$ and $V_2 \supseteq V'$,  
$(\rho(g) V_1) \cap V_2 \supseteq V' \neq 0$ for each $g \in G$. 
In this case, $(\rho(g) V_1) + V_2 \neq V$ for any $g \in G$, and hence $\rho$ can not be $m$-thick. 
If $\ell \ge \max(m, n-m)$, then  let us take subspaces $V_1, V_2$ of $V$ such that 
$V_1 \subseteq V'$ and $V_2 \subseteq V'$. 
Since $(\rho(g)V_1) + V_2 \subseteq V' \neq V$, $\rho$ is not $m$-thick. 
In the case $\min(m, n-m) < \ell < \max(m, n-m)$, we may assume that $n-m \ge m$ because  
$m$-thickness and $(n-m)$-thickness are equivalent. 
Then let us take subspaces 
$V_1, V_2$ of $V$ such that $V_1 \subseteq V' \subseteq V_2$. 
Since $\rho(g)V_1 \subseteq V' \subseteq V_2$, $(\rho(g)V_1) \cap V_2 = \rho(g)V_1 \neq 0$.  
Hence $(\rho(g)V_1) + V_2 \neq V$ for each $g \in G$, which implies $\rho$ is 
not $m$-thick.    
\qed 

\begin{corollary}\label{cor:imply}
For any finite-dimensional 
representation of a group $G$, the following implications hold: 
\begin{eqnarray} 
\mbox{dense} \Rightarrow \mbox{thick} \Rightarrow irreducible. 
\end{eqnarray} 
\end{corollary}
 
\begin{corollary}
Assume that $\dim V \le 3$. 
Then for a representation $\rho : G \to {\rm GL}(V)$, the following 
implications hold:  
\begin{eqnarray} 
\mbox{dense} \Leftrightarrow \mbox{thick} \Leftrightarrow irreducible. 
\end{eqnarray} 
\end{corollary}

{\it Proof}. 
The statement follows from that the three conditions above are 
equivalent to $1$-dense ($1$-thick, or irreducible) when $\dim V \le 3$.  
\hfill $\Box$

\bigskip

Now we consider a criterion for a 
representation to be $m$-thick. 
This criterion of $m$-thickness will be used for describing 
the moduli of absolutely thick representations. 
Before introducing it, we need the following definition. 

\begin{definition}\label{def:realizable}\rm
Let $V$ be a finite-dimensional vector space over a field $k$.  
For a vector subspace $W \subseteq \wedge^m V$, we say that $W$ is 
{\it realizable} over $k$ if there exist 
$v_1, v_2, \ldots, v_m \in V$ such that $0 \neq v_1 \wedge v_2 \wedge 
\cdots \wedge v_m \in W$. 
For an $m$-dimensional subspace $V'$ of $V$ with $0 < m < n$, we can consider a point  
$[ \wedge^m V' ] \in {\mathbb P}_{\ast}(\wedge^m V)$. 
In the sequel, we identify $[ \wedge^m V' ]$ with   
a non-zero vector $\wedge^m V' \in \wedge^m V$ (which is determined 
by $[ \wedge^m V' ]$ up to scalar) for simplicity.  
It is obvious that $W$ is realizable if and only if 
$W$ contains a non-zero vector $\wedge^m V'$ 
obtained by an $m$-dimensional subspace $V'$ of $V$ over $k$.  
%For a vector subspace $W \subseteq \wedge^m V$, we say that $W$ is 
%{\it realizable} over $k$ if $W$ contains a non-zero vector $\wedge^m V'$ 
%obtained by an $m$-dimensional subspace $V'$ of $V$ over $k$. 
\end{definition}

The following proposition gives a criterion of $m$-thickness. 

\begin{proposition}\label{prop:condofm-thick}
Let $\rho : G \to {\rm GL}(V)$ be an $n$-dimensional representation of a group $G$. 
For $0 < m <n$, $\rho$ is not $m$-thick if and only if there exist 
$G$-invariant realizable vector subspaces $W_1 \subset \wedge^m V$ and $W_2 
\subset \wedge^{n-m} V$ 
such that $W_1^{\perp} = W_2$. 
\end{proposition}

{\it Proof}. 
Suppose that $\rho$ is not $m$-thick. Then there exist vector subspaces $V_1, V_2$ of $V$ 
with $\dim V_1 =m$ and $\dim V_2 = n-m$ such that $(\rho(g)V_1) + V_2 \neq V$ for any 
$g \in G$. Let us consider the vector $\wedge^m V_1 \in \wedge^m V$
determined by $V_1$ up to scalar multiplication. 
The condition implies that vectors $\{ (\wedge^m \rho)(g)(\wedge^m V_1) \mid g \in G \}$ 
span a non-trivial $G$-invariant subspace $W_1 \subset \wedge^m V$. 
Of course, $W_1$ is realizable. 
Set $W_2 := W_1^{\perp} \subset \wedge^{n-m}V$. Note that $\wedge^{n-m} V_2 \in W_2$. 
The subspace $W_2$ is a non-trivial $G$-invariant realizable subspace. 
Hence we have proved the ``only if" part. 

Conversely, suppose that there exist 
$G$-invariant realizable vector subspaces $W_1 \subseteq \wedge^m V$ and $W_2 
\subseteq \wedge^{n-m} V$ such that $W_1^{\perp} = W_2$. 
Since $W_1$ and $W_2$ are realizable, there exist an $m$-dimensional subspace
$V_1 \subseteq V$ and an $(n-m)$-dimensional subspace $V_2 \subseteq V$ such 
that $\wedge^m V_1 \in W_1$ and $\wedge^{n-m}V_2 \in W_2$. 
For each $g \in G$, 
the vector $(\wedge^m \rho)(g) (\wedge^m V_1)$ is contained in $W_1$, and  
hence $(\wedge^m\rho)(g) (\wedge^m V_1) \wedge (\wedge^{n-m} V_2) = 0$.  
This implies that $(\rho(g)V_1) + V_2 \neq V$ for each $g \in G$. 
Therefore $\rho$ is not $m$-thick.  
\hfill $\Box$

\begin{remark}\rm
Furthermore, we also see that $\rho$ is not $m$-thick if and only if 
there exist a non-zero 
$G$-invariant realizable subspace $W_1 \subset \wedge^m V$ 
and an $(n-m)$-dimensional subspace $V'$ of $V$ such that 
$\wedge^{n-m} V' \in W_1^{\perp}$. 
%Indeed, we can choose $W_2$ as the $G$-invariant subspace of $\wedge^{n-m} V$ 
%generated by $\wedge^{n-m} V'$.  
\end{remark}

Let us define absolute thickness and absolute denseness. 
We will construct the moduli spaces of 
absolutely thick representations and absolutely 
dense representations in the next section. 

\begin{definition}\rm
Let $G$ be a group. Let $V$ be an $n$-dimensional vector space over 
a field $k$. We say that a representation $\rho : G \to {\rm GL}(V)$ is 
{\it absolutely $m$-thick}  
if $\rho\otimes \overline{k} : 
G \to {\rm GL}(V\otimes \overline{k})$ is $m$-thick, where 
$\overline{k}$ is an algebraic closure of $k$. 
We also say that $\rho$ is 
{\it absolutely thick} if $\rho$ is absolutely $m$-thick for each $0 < m < n$.  
\end{definition}

\begin{definition}\rm
Let $G$ be a group. Let $V$ be an $n$-dimensional vector space over 
a field $k$. We say that a representation $\rho : G \to {\rm GL}(V)$ is 
{\it absolutely $m$-dense} if $\rho\otimes \overline{k} : 
G \to {\rm GL}(V\otimes \overline{k})$ is $m$-dense, where 
$\overline{k}$ is an algebraic closure of $k$. 
We also say that $\rho$ is 
{\it absolutely dense} if $\rho$ is absolutely $m$-dense for each $0 < m < n$.   
\end{definition}

\begin{remark}\label{remark:extfield}\rm 
Let $K$ be an extension field of $k$. 
If $\rho\otimes_{k}K : G \to {\rm GL}(V\otimes_{k}K)$ 
is $m$-thick (resp. $m$-dense), then $\rho$ is also 
$m$-thick (resp. $m$-dense). 
In particular, if $\rho$ is absolutely $m$-thick (resp. absolutely $m$-dense), then 
$\rho$ is $m$-thick (resp. $m$-dense). 
\end{remark}

\begin{proposition}\label{prop:algclosed-dense} 
For an $n$-dimensional group representation $\rho : G \to {\rm GL}(V)$, 
the following conditions are equivalent: 
\begin{enumerate}
\item $\rho$ is absolutely $m$-dense, in other words, 
$(\wedge^m \rho) \otimes_{k}\overline{k} : G \to {\rm GL}(\wedge^m V \otimes_{k} \overline{k})$ is irreducible, where 
$\overline{k}$ is an algebraically closure of $k$. 
\item $(\wedge^m \rho) \otimes_{k}K : G \to {\rm GL}(\wedge^m V \otimes_{k} K)$ is irreducible 
for some algebraically closed field $K$ containing $k$. 
\item $(\wedge^m \rho) \otimes_{k}K : G \to {\rm GL}(\wedge^m V \otimes_{k} K)$ is irreducible 
for any algebraically closed field $K$ containing $k$. 
\end{enumerate} 
\end{proposition} 

{\it Proof}. 
The statement follows from that all conditions above are 
equivalent to the condition that $\wedge^m \rho$ is absolutely 
irreducible. 
\qed 

\bigskip 

In Theorem \ref{th:absm-thick}, 
we will obtain the same result on absolute $m$-thickness 
as Proposition~\ref{prop:algclosed-dense}.  

%%%%%%%%%%%%%%%%%%%%%%%%%%%%%%%%%%%%%%%%%%%%%%%%%%%%%%%%%%%%%%%%%%%%%%%%%%%%%%%%%%%%%%%%%%%%%%%%
%
%   The 3rd Section 
%   The moduli of thick representations 
%
%%%%%%%%%%%%%%%%%%%%%%%%%%%%%%%%%%%%%%%%%%%%%%%%%%%%%%%%%%%%%%%%%%%%%%%%%%%%%%%%%%%%%%%%%%%%%%%%

\section{The moduli of absolutely thick representations}
In this section, we show that absolute thickness is an open condition in 
the representation variety. (For representation varieties,  
see \cite{Nkmt00} )  

Let ${\rm Rep}_n(G)$ be the representation variety of degree $n$ 
for a group $G$ over ${\mathbb Z}$.  
The representation variety represents the following contravariant functor from the category of schemes 
to the category of sets: 
\[
\begin{array}{cccl}
{\rm Rep}_n(G) : & ({\bf Sch})^{op} & \to & ({\bf Sets}) \\  
 & X & \mapsto & 
\left\{ 
%\rho    
%\mid 
 \mbox{ a group representation } 
\rho : G \to  {\rm GL}_n(\Gamma(X, {\mathcal O}_X))  
\right\},   
\end{array} 
\]
where $\Gamma(X, {\mathcal O}_X)$ is the ring of global sections on $X$. 
The representation variety ${\rm Rep}_n(G)$ has the universal $n$-dimensional representation $\tilde{\rho}$ of $G$. 
%on ${\rm Rep}_n(G) \times {\mathbb A}^n_{\mathbb Z}$. 
% 
Let ${\rm Gr}(d, {\mathbb A}^n_{\mathbb Z})$ be the Grassmann scheme over 
${\mathbb Z}$ representing the contravariant functor 
\[
\begin{array}{cccl}
{\rm Gr}(d, {\mathbb A}^n_{\mathbb Z}) : & ({\bf Sch})^{op} & \to & ({\bf Sets}) \\  
 & X & \mapsto & 
\left\{ 
W \;   
\begin{array}{|l}
W \subseteq {\mathcal O}_X^{\oplus n} \mbox{ is a subbundle of rank $d$ } \\ 
\end{array}
\right\}.  
\end{array} 
\]
Let us define a subfunctor $X(d, n; G)$ of ${\rm Rep}_n(G) \times {\rm Gr}(d, {\mathbb A}^n_{\mathbb Z})$ 
for $0 < d < n$ by 
\[
\begin{array}{cccl}
X(d, n; G) : & ({\bf Sch})^{op} & \to & ({\bf Sets}) \\  
 & X & \mapsto & 
\left\{ 
(\rho, W) \;   
\begin{array}{|l}
\rho : G \to  {\rm GL}_n(\Gamma(X, {\mathcal O}_X)), \\
W \subseteq {\mathcal O}_X^{\oplus n} \mbox{ is a subbundle of rank $d$, } \\ 
\mbox{ and } \rho(G) W \subseteq W 
\end{array}
\right\}.  
\end{array} 
\]
We show that $X(d, n; G)$ is a closed subscheme of 
${\rm Rep}_n(G) \times {\rm Gr}(d, {\mathbb A}^n_{\mathbb Z})$. 

\begin{lemma}\label{lemma:3-1}
For $d=1$, $X(d, n; G)$ is a closed subscheme of 
${\rm Rep}_n(G) \times {\rm Gr}(d, {\mathbb A}^n_{\mathbb Z})$. 
\end{lemma}

{\it Proof}. 
The Grassmann scheme ${\rm Gr}(1, {\mathbb A}^n_{\mathbb Z})$ can be  
regarded as ${\mathbb P}_{\ast}({\mathbb A}^n_{\mathbb Z}) := 
\{ [w] \mid w \mbox{ is a non-zero ``vector" of } {\mathbb A}^n_{\mathbb Z} \}$.   
Then 
\begin{eqnarray*}
X(1, n; G) & = & \{ (\rho, [w]) \mid w \mbox{ is a non-zero $\rho(G)$-eigenvector } %of } %{\mathbb A}^n
 \} \\
  & = & \displaystyle \cap_{g \in G} 
\{ (\rho, [w]) \mid w \mbox{ is a non-zero $\rho(g)$-eigenvector } \}.    
\end{eqnarray*}
The condition that $w \in {\mathbb A}^n$ is a $\rho(g)$-eigenvector 
can be written by the equations that all $2$-minors of 
the $n \times 2$ matrix $({\rho} (g)w, w)$ are $0$. 
Hence the subfunctor $X(1, n; G)$ is a closed subscheme of 
${\rm Rep}_n(G) \times {\rm Gr}(1, {\mathbb A}^n_{\mathbb Z})$.
\qed 

\begin{proposition}\label{prop:3-2} 
For $0 < d < n$, $X(d, n; G)$ is a closed subscheme of 
${\rm Rep}_n(G) \times {\rm Gr}(d, {\mathbb A}^n_{\mathbb Z})$. 
\end{proposition}

{\it Proof}. 
The statement is true for $d=1$ by Lemma \ref{lemma:3-1}. 
For $0 < d < n$, 
by taking the exterior, we get the morphism 
\[
\begin{array}{cccc} 
\Phi : &{\rm Rep}_n(G) \times {\rm Gr}(d, {\mathbb A}^n_{\mathbb Z}) & \to 
& {\rm Rep}_{\binom{n}{d}}(G) \times 
{\rm Gr}(1, \wedge^d {\mathbb A}^{n}_{\mathbb Z}) \\
 & (\rho, W) & \mapsto & (\wedge^d \rho, \wedge^d W). 
\end{array}
\]
The subfunctor $X(d, n; G)$ can be obtained by taking the pull-back of the closed subscheme  
$X(1, \binom{n}{d}; G)$ of ${\rm Rep}_{\binom{n}{d}}(G) \times 
{\rm Gr}(1, \wedge^d {\mathbb A}^{n}_{\mathbb Z})$  
by $\Phi$. Hence $X(d, n; G)$ is a closed subscheme of 
${\rm Rep}_n(G) \times {\rm Gr}(d, {\mathbb A}^n_{\mathbb Z})$.  
\hfill $\Box$ 

\bigskip 

Let $0 < m < n$. 
The universal representation $\tilde{\rho}$ on ${\rm Rep}_n(G)$ 
induces an  %family of 
$\binom{n}{m}$-dimensional representation 
$\wedge^m \tilde{\rho}$ on   
${\rm Rep}_{\binom{n}{m}}(G)$. 
This correspondence gives us the canonical morphism 
$\wedge^m : {\rm Rep}_n(G) \to 
{\rm Rep}_{\binom{n}{m}}(G)$  
by $\rho \mapsto \wedge^m \rho$.  
For $0 < d < \binom{n}{m}$, we define
the subfunctor $Y(d, \wedge^m(n); G)$ of 
${\rm Rep}_n(G) \times {\rm Gr}(d, \wedge^m {\mathbb A}^{n}_{\mathbb Z})$ by 
\[
\begin{array}{cccl} 
Y(d, \wedge^m(n); G) : & ({\bf Sch})^{op} & \to & ({\bf Sets}) \\  
 & X & \mapsto & 
\left\{ 
(\rho, W) \;   
\begin{array}{|l}
W \subseteq \wedge^m {\mathcal O}_X^{\oplus n} \mbox{ is a $(\wedge^m \rho)(G)$-invariant } \\ 
\mbox{ subbundle of rank $d$ } \\ 
\end{array}
\right\}.  
\end{array} 
\]
Let us define    
$ \phi := \wedge^m \times id : {\rm Rep}_n(G) \times {\rm Gr}(d, \wedge^m {\mathbb A}^{n}_{\mathbb Z}) \to 
{\rm Rep}_{\binom{n}{m}}(G) \times {\rm Gr}(d, \wedge^m {\mathbb A}^{n}_{\mathbb Z})$ 
by $(\rho, W) \mapsto (\wedge^m \rho, W)$.  
The subfunctor $Y(d, \wedge^m(n); G)$ is obtained by taking 
the pull-back of the closed subscheme $X(d, \binom{n}{m}; G)$ of 
${\rm Rep}_{\binom{n}{m}}(G) \times {\rm Gr}(d, \wedge^m {\mathbb A}^{n}_{\mathbb Z})$ by 
$\phi$. Hence the subfunctor $Y(d, \wedge^m(n); G)$ 
is a closed subscheme of 
${\rm Rep}_n(G) \times {\rm Gr}(d, \wedge^m {\mathbb A}^{n}_{\mathbb Z})$.

We define the subfunctor $Y(d, \wedge^m(n),  \wedge^{n-m}(n); G)$ 
of ${\rm Rep}_n(G) \times {\rm Gr}(d, \wedge^m{\mathbb A}^{n}_{\mathbb Z}) \times 
{\rm Gr}(\binom{n}{m}-d, \wedge^{n-m}{\mathbb A}^{n}_{\mathbb Z})$ by 
\[
\begin{array}{l} 
 Y(d, \wedge^m(n),  \wedge^{n-m}(n); G) :  ({\bf Sch})^{op}  \to  ({\bf Sets}) \\  
\begin{array}{cccl} 
% Y(d, \wedge^m(n),  \wedge^{n-m}(n); G) :  ({\bf Sch})  \to  ({\bf Sets}) \\  
 & X & \mapsto & 
\left\{ 
(\rho, W_1, W_2) \;   
\begin{array}{|l}
W_1 \subseteq \wedge^{m} {\mathcal O}_X^{\oplus n} \mbox{ is a $(\wedge^m \rho)(G)$-invariant } \\ 
\mbox{ subbundle of rank $d$, and} \\ 
W_2 \subseteq \wedge^{n-m} {\mathcal O}_X^{\oplus n} \mbox{ is a $(\wedge^{n-m} \rho)(G)$-invariant } \\ 
\mbox{ subbundle of rank $\binom{n}{m} - d$ } \\ 
\end{array}
\right\}.  
\end{array} 
\end{array} 
\]
Set 
$X_{n, m, d}(G) : = {\rm Rep}_n(G) \times 
{\rm Gr}(d, \wedge^{m} {\mathbb A}^{n}_{\mathbb Z}) \times 
{\rm Gr}(\binom{n}{m}-d, \wedge^{n-m} {\mathbb A}^{n}_{\mathbb Z})$. Let us consider the two projections 
\[
\begin{array}{ccl}
\phi_1 & : & 
X_{n, m, d}(G) \to 
{\rm Rep}_n(G) \times {\rm Gr}(d, \wedge^m {\mathbb A}^{n}_{\mathbb Z}) \\
\phi_2 & : &  X_{n, m, d}(G) \to 
{\rm Rep}_n(G) \times {\rm Gr}(\binom{n}{m}-d, \wedge^{n-m} 
{\mathbb A}^{n}_{\mathbb Z}).  
\end{array}
\] 
Take the pull-backs $\phi_1^{-1}(Y(d, \wedge^m(n); G))$ and 
$\phi_2^{-1}(Y(\binom{n}{m}-d, \wedge^{n-m}(n); G))$.   
The subfunctor $Y(d, \wedge^m(n),  \wedge^{n-m}(n); G)$  
can be obtained as the intersection of these two pull-backs. 
Therefore $Y(d, \wedge^m(n),  \wedge^{n-m}(n); G)$ is 
a closed subscheme of $X_{n, m, d}(G)$.

Set ${\rm Gr}_{n, m, d} := {\rm Gr}(d, \wedge^{m} {\mathbb A}^{n}_{\mathbb Z}) 
\times {\rm Gr}(\binom{n}{m}-d, \wedge^{n-m} {\mathbb A}^{n}_{\mathbb Z})$.  
Let us consider the perfect pairing on ${\rm Gr}_{n, m, d}$:  
\[ 
\langle \; , \; \rangle : (\wedge^{m} {\mathcal O}_{{\rm Gr}_{n, m, d}}^{\oplus n})  
\otimes_{{\mathcal O}_{{\rm Gr}_{n, m, d}}} 
(\wedge^{n-m} {\mathcal O}_{{\rm Gr}_{n, m, d}}^{\oplus n})  
\to \wedge^{n} {\mathcal O}_{{\rm Gr}_{n, m, d}}^{\oplus n} 
\cong {\mathcal O}_{{\rm Gr}_{n, m, d}}  
\]
defined by $\langle x, y \rangle := x \wedge y$. 
We define the subfunctor ${\rm Gr}_{n, m, d}^{\perp}$ of ${\rm Gr}_{n, m, d}$ by 
\[
{\rm Gr}_{n, m, d}^{\perp} := \{ (W_1, W_2) \in {\rm Gr}_{n, m, d} \mid 
W_1^{\perp} = W_2 \}. 
\] 
For each point $p = (W_1, W_2) \in {\rm Gr}_{n, m, d}$, choose a neighbourhood $U$ of $p$ and 
sections $\{ e_i \}, \{ f_j \}$ on $U$ 
such that $\langle e_1, e_2, \ldots, e_d \rangle$ is the universal subbundle of 
$\wedge^m {\mathcal O}_{{\rm Gr}_{n, m, d}}^{\oplus n}$ of rank $d$ on $U$ 
and $W_2 = \langle f_1, f_2, \ldots, f_{\binom{n}{m}-d} \rangle$ is 
the universal subbundle of $\wedge^{n-m} {\mathcal O}_{{\rm Gr}_{n, m, d}}^{\oplus n}$ 
of rank $\binom{n}{m}-d$ on $U$.  
The equations $\langle e_i, f_j \rangle = 0$ define a closed subscheme structure on 
${\rm Gr}_{n, m, d}^{\perp}$. 
Hence ${\rm Gr}_{n, m, d}^{\perp}$ is a closed subscheme of 
${\rm Gr}_{n, m, d}$.  

Let us denote by $\phi_3 : X_{n, m, d}(G) \to {\rm Gr}_{n, m, d}$ the 
canonical projection. Taking the intersection of $Y(d, \wedge^m(n),  \wedge^{n-m}(n); G)$
with the pull-back $\phi_3^{-1}({\rm Gr}_{n, m, d}^{\perp})$, we obtain a  
closed subscheme $Y(d, \wedge^m(n),  \wedge^{n-m}(n); G)^{\perp}$ of 
$Y(d, \wedge^m(n),  \wedge^{n-m}(n); G)$. 
The closed subscheme $Y(d, \wedge^m(n),  \wedge^{n-m}(n); G)^{\perp}$
represents the contravariant functor 
\[
\begin{array}{cccl} 
%Y(d, \wedge^m(n), \wedge^{n-m}(n); G) : 
& ({\bf Sch})^{op} & \to & ({\bf Sets}) \\  
 & X & \mapsto & 
\left\{ 
(\rho, W_1, W_2) \in Y(d, \wedge^m(n), \wedge^{n-m}(n); G)(X) \;   
\begin{array}{|l}
W_1^{\perp} = W_2 \\ 
\end{array}
\right\}.  
\end{array} 
\]

\bigskip 

For proving openness of absolute $m$-thickness, 
we show that realizable subspaces form a closed subset in the Grassmann scheme.  
We set 
\[
{\rm Gr}(d, \wedge^m {\mathbb A}^{n}_{\mathbb Z})_{real} := 
\{ W \in {\rm Gr}(d, \wedge^m {\mathbb A}^{n}_{\mathbb Z}) \mid 
W \mbox{ is realizable } \}. 
\]  
More precisely, for a point $x \in {\rm Gr}(d, \wedge^m{\mathbb A}^{n}_{\mathbb Z})$, 
$x \in {\rm Gr}(d, \wedge^m {\mathbb A}^{n}_{\mathbb Z})_{real}$ if 
and only if there exists an extension field $K$ of the residue field $k(x)$ of $x$ 
such that the $d$-dimensional subspace $W \subseteq \wedge^m K^{n}$ associated to $x$ is realizable over $K$. 

\begin{prop}\label{prop:grassrealclosed}
The subset ${\rm Gr}(d, \wedge^m {\mathbb A}^{n}_{\mathbb Z})_{real}$ can be regarded as a closed subscheme of ${\rm Gr}(d, \wedge^m {\mathbb A}^{n}_{\mathbb Z})$.  
\end{prop} 

{\it Proof}.
Let us consider the closed subscheme 
\[
  {\rm Flag}(1, d, \wedge^m {\mathbb A}^{n}_{\mathbb Z}) := \{ ([v], W) \in 
{\mathbb P}_{\ast}(\wedge^{m} {\mathbb A}^{n}_{\mathbb Z}) \times 
{\rm Gr}(d, \wedge^{m} {\mathbb A}^{n}_{\mathbb Z})  
 \mid v \in W \} 
\]
of ${\mathbb P}_{\ast}(\wedge^{m} {\mathbb A}^{n}_{\mathbb Z}) \times 
{\rm Gr}(d, \wedge^{m} {\mathbb A}^{n}_{\mathbb Z}) = 
{\rm Gr}(1, \wedge^{m} {\mathbb A}^{n}_{\mathbb Z}) \times 
{\rm Gr}(d, \wedge^{m} {\mathbb A}^{n}_{\mathbb Z})$.  
The scheme ${\mathbb P}_{\ast}(\wedge^{m} {\mathbb A}^{n}_{\mathbb Z})$ has a closed subscheme 
${\rm Gr}(m, {\mathbb A}^n_{\mathbb Z})$. Then we obtain the pull-back 
$p_{1}^{-1}({\rm Gr}(m, {\mathbb A}^n_{\mathbb Z}))$ of ${\rm Gr}(m, {\mathbb A}^n_{\mathbb Z})$ 
by the first projection $p_1 : {\rm Flag}(1, d, \wedge^{m} {\mathbb A}^{n}_{\mathbb Z}) 
\to {\mathbb P}_{\ast}(\wedge^{m} {\mathbb A}^{n}_{\mathbb Z})$.  
The subset ${\rm Gr}(d, \wedge^{m} {\mathbb A}^{n}_{\mathbb Z})_{real}$ is 
the image $p_2( p_{1}^{-1}({\rm Gr}(m, {\mathbb A}^n_{\mathbb Z})) )$ of the closed subscheme 
$p_{1}^{-1}({\rm Gr}(m, {\mathbb A}^n_{\mathbb Z}))$ by the second projection 
$p_2 : {\rm Flag}(1, d, \wedge^{m} {\mathbb A}^{n}_{\mathbb Z}) \to 
{\rm Gr}(d, \wedge^{m} {\mathbb A}^{n}_{\mathbb Z})$.
The projection $p_2$ is proper, and hence we can define a closed subscheme structure 
on ${\rm Gr}(d, \wedge^{m} {\mathbb A}^{n}_{\mathbb Z})_{real}$.    
\hfill $\Box$

\begin{remark}\label{remark:radon}\rm 
In the proof of Proposition \ref{prop:grassrealclosed},  we also see that   
\begin{eqnarray} 
{\rm Gr}(d, \wedge^m {\mathbb A}^{n}_{\mathbb Z})_{real}
& = & p_2( p_{1}^{-1}({\rm Gr}(m, {\mathbb A}^n_{\mathbb Z})) ). 
\end{eqnarray} 
\end{remark}

\bigskip 

The following proposition gives a characterization of 
${\rm Gr}(d, \wedge^{m} {\mathbb A}^{n}_{\mathbb Z})_{real}$. 

\begin{prop}\label{prop:grassreal} 
Let $x \in {\rm Gr}(d, \wedge^m {\mathbb A}^{n}_{\mathbb Z})$.  
Let $\overline{k(x)}$ be an algebraic closure of 
the residue field $k(x)$ of $x$. 
Then $x \in {\rm Gr}(d, \wedge^m {\mathbb A}^{n}_{\mathbb Z})_{real}$ if and only if 
the corresponding $d$-dimensional subspace $W\otimes_{k(x)} \overline{k(x)} 
\subseteq \wedge^m \overline{k(x)}^{\:n}$ to $x$ is realizable over 
$\overline{k(x)}$. 
\end{prop} 

{\it Proof}. 
Let $x \in {\rm Gr}(d, \wedge^m {\mathbb A}^{n}_{\mathbb Z})$. 
If the corresponding $d$-dimensional subspace $W\otimes_{k(x)} \overline{k(x)} 
\subseteq \wedge^m \overline{k(x)}^{\:n}$ is realizable over 
$\overline{k(x)}$, then 
there exists a $\overline{k(x)}$-rational point 
of $p_1^{-1}({\rm Gr}(m, {\mathbb A}^n_{\mathbb Z})) \subseteq {\rm Flag}(1, d, \wedge^m {\mathbb A}^{n}_{\mathbb Z})$ 
whose image by $p_2$ corresponds to $W\otimes_{k(x)} \overline{k(x)}$. 
Then $x \in p_2(p_1^{-1}({\rm Gr}(m, {\mathbb A}^n_{\mathbb Z}))) = 
{\rm Gr}(d, \wedge^m {\mathbb A}^{n}_{\mathbb Z})_{real}$. 

Conversely, suppose that $x \in {\rm Gr}(d, \wedge^m {\mathbb A}^{n}_{\mathbb Z})_{real}$.  
Setting $\phi := p_2\!\!\mid_{p_1^{-1}({\rm Gr}(m, {\mathbb A}^n_{\mathbb Z}))}$, we have 
the following commutative diagram which is a fibre product: 
\[
\begin{array}{ccc}
p_1^{-1}({\rm Gr}(m, {\mathbb A}^n_{\mathbb Z})) & 
\stackrel{\phi}{\to} & {\rm Gr}(d, \wedge^m {\mathbb A}^{n}_{\mathbb Z})_{real} \\
\uparrow & & \uparrow \\
\phi^{-1}(x) & \to & {\rm Spec}\:k(x)  . 
\end{array}
\]
Since $\phi$ is of finite type, so is $\phi^{-1}(x) \to {\rm Spec}\:k(x)$.  
Note that $\phi^{-1}(x) \neq \emptyset$ 
by the definition of ${\rm Gr}(d, \wedge^m {\mathbb A}^{n}_{\mathbb Z})_{real}$.  
Then there exists a $\overline{k(x)}$-rational point 
of $\phi^{-1}(x)$. 
This implies that the corresponding $d$-dimensional subspace $W\otimes_{k(x)} \overline{k(x)} 
\subseteq \wedge^m \overline{k(x)}^{\:n}$ is realizable over 
$\overline{k(x)}$. 
\qed

\bigskip 

Let $q_2 : X_{n, m, d}(G) \to {\rm Gr}(d, \wedge^{m} {\mathbb A}^n_{\mathbb Z})$ 
and 
$q_3 : X_{n, m, d}(G) \to {\rm Gr}(\binom{n}{m}-d, \wedge^{n-m} {\mathbb A}^n_{\mathbb Z})$ 
be the second and the third projections. 
We denote by $Y(d, \wedge^m(n),  \wedge^{n-m}(n); G)^{\perp}_{real}$ the intersection of  
$Y(d, \wedge^m(n),  \wedge^{n-m}(n); G)^{\perp}$ with  
$q_{2}^{-1}({\rm Gr}(d, \wedge^{m} {\mathbb A}^{n}_{\mathbb Z})_{real}) \cap 
q_{3}^{-1}({\rm Gr}(\binom{n}{m}-d, \wedge^{n-m} {\mathbb A}^{n}_{\mathbb Z})_{real})$. 
By Proposition \ref{prop:grassrealclosed}, $Y(d, \wedge^m(n),  \wedge^{n-m}(n); G)^{\perp}_{real}$ 
can be regarded as a closed subscheme of $Y(d, \wedge^m(n),  \wedge^{n-m}(n); G)^{\perp}$. 

%\bigskip 

\begin{prop}\label{prop:y}
Let $x =(\rho, W_1, W_2) \in Y(d, \wedge^m(n),  \wedge^{n-m}(n); G)^{\perp}$. 
Let $\overline{k(x)}$ be an algebraic closure of 
the residue field $k(x)$ of $x$. Then   
$x \in Y(d, \wedge^m(n),  \wedge^{n-m}(n); G)^{\perp}_{real}$ if and only if 
$W_1\otimes_{k(x)}\overline{k(x)}$ and 
$W_2\otimes_{k(x)}\overline{k(x)}$ are realizable over $\overline{k(x)}$. 
\end{prop} 

{\it Proof}.
By the definition, $x \in Y(d, \wedge^m(n),  \wedge^{n-m}(n); G)^{\perp}_{real}$ 
if and only if $q_{2}(x) \in {\rm Gr}(d, \wedge^{m} {\mathbb A}^{n}_{\mathbb Z})_{real}$ and 
$q_{3}(x) \in {\rm Gr}(\binom{n}{m}-d, \wedge^{n-m} {\mathbb A}^{n}_{\mathbb Z})_{real}$. 
It follows from Proposition \ref{prop:grassreal} that   
this condition is equivalent to that  
$W_1\otimes_{k(x)}\overline{k(x)}$ and 
$W_2\otimes_{k(x)}\overline{k(x)}$ are realizable over $\overline{k(x)}$. 
\qed

\bigskip 

Let $q_1 : X_{n, m, d}(G) \to {\rm Rep}_n(G)$ be the first projection.  
Since $q_1$ is proper and $Y(d, \wedge^m(n),  \wedge^{n-m}(n); G)^{\perp}_{real}$ 
is a closed subscheme of $X_{n, m, d}(G)$, ${\rm Rep}_n(G)$ has a 
closed subscheme 
$q_1(Y(d, \wedge^m(n),  \wedge^{n-m}(n); G)^{\perp}_{real})$.  

\begin{prop}\label{prop:q1}
Let $x \in {\rm Rep}_n(G)$. 
Then $x \in q_1(Y(d, \wedge^m(n),  \wedge^{n-m}(n); G)^{\perp}_{real})$ if and only if 
there exist $G$-invariant realizable subspaces $W_1 
\subseteq \wedge^{m} \overline{k(x)}^{\:n}$ and 
$W_2 \subseteq \wedge^{n-m} \overline{k(x)}^{\:n}$ 
with respect to the corresponding representation $\rho_{x}\otimes_{k(x)} \overline{k(x)} 
: G \to {\rm GL}_n(\overline{k(x)})$ 
such that $\dim W_1 = d , \dim W_2 = \binom{n}{m}-d$, and $W_1^{\perp} = W_2$. 
\end{prop} 

{\it Proof}.
First, we prove the ``if" part. 
Suppose that there exist such $W_1$ and $W_2$. 
Then we have a $\overline{k(x)}$-rational point of 
$Y(d, \wedge^m(n),  \wedge^{n-m}(n); G)^{\perp}_{real}$ whose 
image by $q_1$ corresponds to $x$. 
Hence $x \in q_1(Y(d, \wedge^m(n),  \wedge^{n-m}(n); G)^{\perp}_{real})$. 

Next, we prove the ``only if" part. Let 
$x \in q_1(Y(d, \wedge^m(n),  \wedge^{n-m}(n); G)^{\perp}_{real})$. 
Set $\psi = q_1\!\!\mid_{Y(d, \wedge^m(n),  \wedge^{n-m}(n); G)^{\perp}_{real}} 
: Y(d, \wedge^m(n),  \wedge^{n-m}(n); G)^{\perp}_{real} \to {\rm Rep}_n(G)$. 
Since $\psi$ is of finite type, so is $\psi^{-1}(x) \to 
{\rm Spec}\: k(x)$.   
The fibre $\psi^{-1}(x)$ is not empty, and hence there exist 
$W_1$ and $W_2$ with the desired property by Proposition \ref{prop:y}. 
\qed

\bigskip 

We can prove the following theorem on absolute $m$-thickness 
by Proposition \ref{prop:q1}.

\begin{theorem}\label{th:absm-thick} 
Let $\rho : G \to {\rm GL}_n(k)$ be an $n$-dimensional representation of $G$ over a field $k$. 
For $0 < m < n$, the following conditions are equivalent:   
\begin{enumerate}
\item\label{cond:1} $\rho$ is absolutely $m$-thick, in other words, $\rho \otimes_{k} \overline{k}$ is 
$m$-thick for an algebraic closure $\overline{k}$ of $k$.   
\item\label{cond:2} $\rho \otimes_{k} K$ is $m$-thick for some algebraically closed field $K$ over $k$. 
\item\label{cond:3} $\rho \otimes_{k} K$ is $m$-thick for any algebraically closed field $K$ over $k$. 
\end{enumerate}
\end{theorem} 

{\it Proof}.
It is obvious that (\ref{cond:3}) $\Rightarrow$ (\ref{cond:1}) 
and that (\ref{cond:1}) $\Rightarrow$ (\ref{cond:2}). 
Let us show that (\ref{cond:2}) $\Rightarrow$ (\ref{cond:3}). 
Assume that $\rho \otimes_{k} K$ is $m$-thick for some 
algebraically closed field $K$ over $k$. 
Note that $\rho \otimes_{k}\overline{k}$ is also $m$-thick by Remark~\ref{remark:extfield}. 
Suppose that $\rho \otimes_{k} K'$ is not $m$-thick for some 
algebraically closed field $K'$ over $k$. 
Let $x$ be the $k$-rational point of ${\rm Rep}_n(G)$ associated to $\rho$. 
By Proposition \ref{prop:condofm-thick},  
there exists a $K'$-rational point of 
$Y(d, \wedge^m(n),  \wedge^{n-m}(n); G)^{\perp}_{real}$ for some $d$ 
whose image by $q_1$ corresponds to $\rho \otimes_{k} K'$. 
Hence $x \in q_1(Y(d, \wedge^m(n),  \wedge^{n-m}(n); G)^{\perp}_{real})$. 
%For any algebraically closed field $K'$ containing $k$, the $K'$-rational 
%point of ${\rm Rep}_n(G)$ associated to $\rho\otimes_{k} K'$ is mapped to 
%a $K'$-rational point of 
%$q_1(Y(d, \wedge^m(n),  \wedge^{n-m}(n); G)^{\perp}_{real})$ 
%by $q_1$.  
Then $\rho \otimes_{k} \overline{k}$ 
is not $m$-thick by Proposition \ref{prop:q1}, which is a contradiction.   
Hence $\rho\otimes_{k} K'$ is $m$-thick for any algebraically closed field $K'$ over $k$.  
Therefore we have shown that 
(\ref{cond:2}) $\Rightarrow$ (\ref{cond:3}). 
\qed 

\bigskip 

Now we show the openness of absolute $m$-thickness. 

\begin{theorem}\label{th:mainth}
Let ${\rm Rep}_n(G)$ be the representation variety of degree $n$ 
for a group $G$ over ${\mathbb Z}$.  
For $0 < m < n$, the absolutely $m$-thick representations in ${\rm Rep}_n(G)$ form an 
open subscheme of ${\rm Rep}_n(G)$. 
In particular, the absolutely thick representations in ${\rm Rep}_n(G)$ form an 
open subscheme of ${\rm Rep}_n(G)$. 
\end{theorem}

{\it Proof}. 
The absolutely $m$-thick representations form the complement of 
\[
\displaystyle 
\bigcup_{0 < d < \binom{n}{m}} q_1(Y(d, \wedge^m(n),  \wedge^{n-m}(n); G)^{\perp}_{real})
\]
in ${\rm Rep}_n(G)$ by Propositions \ref{prop:condofm-thick} and \ref{prop:q1}. 
Since $q_1(Y(d, \wedge^m(n),  \wedge^{n-m}(n); G)^{\perp}_{real})$ is closed for 
each $d$, we can verify the openness of absolute $m$-thickness. 
We can also prove the openness of absolute thickness by considering all $m$.  
\hfill $\Box$

\bigskip 

Let $\repmthick$ be the open subscheme consisting of 
absolutely $m$-thick representations of ${\rm Rep}_n(G)$. 
Let $\repthick$ be the open subscheme consisting of 
absolutely thick representations of ${\rm Rep}_n(G)$. 
The open subschemes $\repmthick$ and 
$\repthick$ 
are contained in the representation variety of 
absolutely irreducible representations 
${\rm Rep}_n(G)_{\rm air}$. We have group actions of the group scheme 
${\rm PGL}_n$ on these schemes  
by the conjugation $\rho \mapsto P^{-1}\rho P$. 
By \cite[Theorem 1.3]{Nkmt00}, there exists a universal geometric quotient 
${\rm Ch}_n(G)_{\rm air}$ of ${\rm Rep}_n(G)_{\rm air}$ by 
${\rm PGL}_n$ and the quotient morphism 
${\rm Rep}_n(G)_{\rm air} \to {\rm Ch}_n(G)_{\rm air}$ 
is a ${\rm PGL}_n$-principal fibre bundle. 
Hence we have the following theorem: 

\begin{theorem}\label{th:modulithick} 
For each $0 < m < n$, there exists a universal geometric quotient 
$\chmthick$ of $\repmthick$ 
by ${\rm PGL}_n$. Moreover, 
there exists a universal geometric quotient 
$\chthick$ of $\repthick$ 
by ${\rm PGL}_n$. 
The quotient morphisms $\repmthick  
\to \chmthick$ and $\repthick 
\to \chthick$ are ${\rm PGL}_n$-principal fibre bundles. 
\end{theorem}

\bigskip 

We also have the same results on absolutely dense representations 
as absolutely thick representations. 

\begin{proposition}
For $0 < m < n$, the absolutely $m$-dense representations in ${\rm Rep}_n(G)$ form an 
open subscheme of ${\rm Rep}_n(G)$. 
In particular, the absolutely dense representations in ${\rm Rep}_n(G)$ form an 
open subscheme of ${\rm Rep}_n(G)$. 
\end{proposition}

{\it Proof}. 
We define the morphism 
$\wedge^m : {\rm Rep}_n(G) \to {\rm Rep}_{\binom{n}{m}}(G)$ 
by $\rho \mapsto \wedge^m \rho$. 
The inverse image of the open 
subscheme ${\rm Rep}_{\binom{n}{m}}(G)_{\rm air}$ by 
$\wedge^m$ coincides with the absolutely 
$m$-dense representations in ${\rm Rep}_n(G)$. 
Hence it is open. 
Considering all $m$, we see that 
the absolutely dense representations in ${\rm Rep}_n(G)$ is also open.
\qed

\bigskip 

Let $\repmdense$ be the open subscheme consisting of 
absolutely $m$-dense representations of ${\rm Rep}_n(G)$. 
Let $\repdense$ be the open subscheme consisting of 
absolutely dense representations of ${\rm Rep}_n(G)$. 
The open subschemes $\repmdense$ and 
$\repdense$ 
are contained in the representation variety of 
absolutely irreducible representations 
${\rm Rep}_n(G)_{\rm air}$. 
In the same way as absolutely thick representations, we have 
the following theorem: 

\begin{theorem}\label{th:modulidense} 
For each $0 < m < n$, there exists a universal geometric quotient 
$\chmdense$ of $\repmdense$ 
by ${\rm PGL}_n$. Moreover, 
there exists a universal geometric quotient 
$\chdense$ of $\repdense$ 
by ${\rm PGL}_n$. The quotient morphisms $\repmdense  
\to \chmdense$ and $\repdense  
\to \chdense$ are ${\rm PGL}_n$-principal fibre bundles. 
\end{theorem}

\bigskip 

Summarizing the results above, we have the following diagrams:  
\[
\begin{array}{ccccc}
\repmdense & \subseteq & 
\repmthick & \subseteq & {\rm Rep}_n(G)_{\rm air} \\
\downarrow & & \downarrow & & \downarrow \\
\chmdense & \subseteq & 
\chmthick & \subseteq & {\rm Ch}_n(G)_{\rm air} \\
\end{array}
\]
and 
\[
\begin{array}{ccccc}
\repdense & \subseteq & 
\repthick & \subseteq & {\rm Rep}_n(G)_{\rm air} \\
\downarrow & & \downarrow & & \downarrow \\
\chdense & \subseteq & 
\chthick & \subseteq & {\rm Ch}_n(G)_{\rm air} . \\
\end{array}
\]

\begin{remark}\rm
For a representation $\rho : G \to {\rm GL}_n(\Gamma(X, {\mathcal O}_X))$ of a group 
$G$ on a scheme $X$, $\rho$ is called {\it absolutely $m$-thick} 
(resp. {\it absolutely thick}) 
if the induced representation $\rho\otimes k(x) : G \to {\rm GL}_n(k(x))$ 
is absolutely $m$-thick (resp. absolutely thick) 
for each $x \in X$, where $k(x)$ is the 
residue field of $x$.  Similarly, $\rho$ is called 
{\it absolutely $m$-dense} (resp. {\it absolutely dense})   
if the induced representation $\rho\otimes k(x) : G \to {\rm GL}_n(k(x))$ 
is absolutely $m$-dense (resp. absolutely dense) 
for each $x \in X$.  
The scheme $\repmthick$ (resp. $\repthick$, $\repmdense$, $\repdense$) represents  
the contravariant functor from the category of schemes to the category of sets 
which maps each scheme to the set of $n$-dimensional 
absolutely $m$-thick (resp. 
absolutely thick, absolutely $m$-dense, absolutely 
dense) representations of $G$ on $X$.      
\end{remark} 

%%%%%%%%%%%%%%%%%%%%%%%%%%%%%%%%%%%%%%%%%%%%%%%%%%%%%%%%%%%%%%%%%%%%%%%%%%%%%%%%%%%%%%%%%%%%%%%%%%%%%%%%%%

\section{Realizable subspaces}

In this section, we discuss realizable subspaces in detail. 
We introduce the $r$-number $r({\wedge^m (n)})$ which 
is closely related to thickness. In some cases, 
we can calculate $r({\wedge^m (n)})$.

\begin{lemma}\label{lemma:grass3}
Let $V$ be an $n$-dimensional vector space over an algebraically closed field $k$.   
Let $W$ be a vector subspace of $\wedge^m V$ with $0 < m < n$.  
If ${\rm codim}\: W \le m(n-m)$, then $W$ is realizable, in other words, there exists an $m$-dimensional vector subspace $V_1$ 
of $V$ such that $\wedge^m V_1 \in W$. 
\end{lemma}

{\it Proof}. 
Remark that $\wedge^m V_1 \in \wedge^m V$ 
can be defined up to scalar multiplication. 
The Grassmann variety ${\rm Gr}(m, V) \subset {\mathbb P}_{\ast}(\wedge^m V)$ has 
dimension $m(n-m)$.   
Since the subspace ${\mathbb P}_{\ast}(W) \subset {\mathbb P}_{\ast}(\wedge^m V)$ 
has codimension $\le m(n-m)$, the intersection ${\mathbb P}_{\ast}(W) 
\cap {\rm Gr}(m, V)$ is not empty. Hence there exists an $m$-dimensional 
subspace $V_1$ such that $\wedge^m V_1 \in W$. 
\qed 

\begin{proposition}\label{prop:condnothick}
Let $V$ be an $n$-dimensional vector space over an algebraically closed field $k$. 
Let $\rho : G \to {\rm GL}(V)$ be a 
representation of a group $G$.   
If $\wedge^m V$ has a $(\wedge^m \rho)(G)$-invariant 
realizable subspace $W$ of $\dim W \le m(n-m)$, then 
$\rho$ is not $m$-thick. 
\end{proposition}

{\it Proof}. 
Let us consider $W^{\perp} \subseteq \wedge^{n-m} V$. 
Since $\dim W \le m(n-m)$, ${\rm codim} W^{\perp} \le m(n-m)$. 
By Lemma \ref{lemma:grass3}, $W^{\perp}$ is realizable. 
Hence $\rho$ is not $m$-thick because of Proposition \ref{prop:condofm-thick}. 
\qed 

\begin{definition}\rm 
For $0 < m < n$, 
 we define the {\it $r$-number} $r({\wedge^m (n)})$ by 
\[
r({\wedge^m (n)}) := \min 
\left\{ \dim W \; 
\begin{array}{|l}
 \mbox{ there exists  
an $n$-dimensional irreducible } \\ 
 \mbox{ representation $\rho : G \to {\rm GL}(V)$ of a group $G$ } \\
 \mbox{ over a field $k$ such that $W$ is a $G$-invariant } \\
\mbox{ realizable subspace of } \wedge^m V  \\
\end{array}
\right\}. 
\] 
For convenience, we set $r({\wedge^0 (n)}) = 1$ and 
$r({\wedge^n (n)}) = 1$ for each positive integer $n$.  
\end{definition}

For a real number $x$, we denote by $[x]$ the largest integer which is equal to or less than $x$.

\begin{proposition}\label{prop:rinequality}
For $0< m < n$, $r({\wedge^m (n)}) \ge [\frac{n-1}{m}] + 1$. 
\end{proposition}

{\it Proof}.
Let $\rho : G \to {\rm GL}(V)$ be an $n$-dimensional irreducible 
representation of a group $G$. 
Let $W \subseteq \wedge^m V$ be a $G$-invariant realizable subspace. 
We show that $\dim W \ge [\frac{n-1}{m}] + 1$. 
Since $W$ is realizable, there exists a basis $e_1, e_2, \ldots, e_n$ of $V$ 
such that $x := e_1 \wedge e_2 \wedge \cdots \wedge e_m \in W$. 
We define $g_i \in G$ for $1 \le i \le [\frac{n-1}{m}]+1$ in the following way: 
Let $g_1 := e \in G$. If $g_i \in G$ is determined for $i \le k$, 
choose $g_{k+1} \in G$ such that $\rho(g_{k+1})e_1$ is not contained in 
the subspace $V_k$ spanned by 
$\{  \rho(g_i)e_j \mid 1 \le i \le k, 1 \le j \le m \}$ of $V$.  
This procedure is possible, since $\dim V_k \le km \le [\frac{n-1}{m}]m < n$ and the 
set $\{ \rho(g)e_1 \mid g \in G \}$ spans $V$ because of the irreducibility of $\rho$. 
In this way, $g_1, g_2, \ldots, g_{[\frac{n-1}{m}]+1}$ can be chosen.   

We claim that 
$(\wedge^m \rho)(g_1)x, (\wedge^m \rho)(g_2)x,
\ldots, (\wedge^m \rho)(g_{[\frac{n-1}{m}]+1})x \in W$ are linearly independent.  
Let $\sum a_i (\wedge^m \rho)(g_{i})x = 0$ for $a_i \in k$.  
By using 
\begin{eqnarray*}
 (\wedge^m \rho)(g_{i})x \wedge \rho(g_{i+1}) e_1 \wedge \rho(g_{i+2}) e_1 \wedge 
\cdots \wedge \rho(g_{[\frac{n-1}{m}]+1}) e_1 \neq 0 \mbox{ and } \\
 (\wedge^m \rho)(g_{i})x \wedge \rho(g_{i})e_1 \wedge \rho(g_{i+1}) e_1 \wedge \rho(g_{i+2}) e_1 \wedge 
\cdots \wedge \rho(g_{[\frac{n-1}{m}]+1}) e_1 = 0,  
\end{eqnarray*} 
we see that $a_i = 0$ for each $i$, which implies the claim.  
Hence $\dim W \ge [\frac{n-1}{m}]+1$.         
\qed 

\begin{corollary}\label{cor:rge2} 
If $0 < m < n$, then $r(\wedge^m(n)) \ge 2$. 
In particular, if $\rho$ is an $n$-dimensional irreducible 
representation, then $\wedge^m \rho$ has no $1$-dimensional 
$G$-invariant realizable subspace. 
\end{corollary}

{\it Proof}. 
The statement follows from that 
$r(\wedge^m(n)) \ge [\frac{n-1}{m}] + 1 \ge 2$. 
\qed 

\bigskip 

%  The following conjecture is not true because we can verify the $n=10, m=5$ case. 
%
%\begin{conjecture}
%Let $\rho :  G \to {\rm GL}(V)$ be an irreducible $n$-dimensional dimensional representation of a group $G$. 
%If the induced representation $\wedge^m V$ has a 
%$G$-invariant realizable subspace $W$ for some $0 < m < n$, then 
%$\dim W \ge {n \choose m} - m(n-m)$.  
%In other words, $r(\wedge^m(n)) \ge {n \choose m} -m(n-m)$. 
%\end{conjecture}
%

If $m$ divides $n$, then we can prove that $r(\wedge^m(n))=\frac{n}{m}$. 
For proving this, we need to make some preparations.

\begin{lemma}\label{lemma:linearalg} 
Let $f : V \to V$ be a linear endomorphism on an $n$-dimensional 
vector space $V$ over a field $k$. Suppose that $f$ has $n$ distinct 
eigenvalues $\alpha_1, \ldots, \alpha_n \in k$. 
Let $e_1, \ldots, e_n \in V$ be eigenvectors 
associated to $\alpha_1, \ldots, \alpha_n$, respectively. 
Then for any $f$-invariant subspace $W$ of $V$, 
there exists a subset $I$ of $\{ 1, 2, \ldots, n \}$ 
such that $W = \oplus_{i \in I} k\cdot e_i$. 
\end{lemma}

{\it Proof}.
For an $f$-invariant subspace $W$, we define a subset $I$ of 
$\{ 1, 2, \ldots, n \}$ by $I := \{ i \mid \mbox{there exists } \sum_{j=1}^n a_j e_j \in W 
\mbox{ such that } a_i \neq 0 \}$. It is clear that $W \subseteq 
\oplus_{i \in I} k\cdot e_i$. 
We show that $W \supseteq 
\oplus_{i \in I} k\cdot e_i$. 
For each $i \in I$, there exists a vector $x = \sum_{j=1}^n a_j e_j 
\in W$ such that $a_i \neq 0$. 
Set $J:= \{ j \mid a_j \neq 0\} = \{ j_1, j_2, \ldots, j_m \}$ and 
$m := \sharp J$. Note that $i \in J$. 
Since $f(x) = \sum_{j \in J} \alpha_j a_j e_j, 
f^2(x) = \sum_{j \in J} \alpha_j^2 a_j e_j, \ldots, f^{m-1}(x) = 
\sum_{j \in J} \alpha_j^{m-1} a_j e_j$, 
we have 
\[
\left(
\begin{array}{c}
x \\
f(x) \\
f^2(x) \\
\vdots \\
f^{m-1}(x) \\
\end{array}
\right) = 
\left(
\begin{array}{cccc}
1 & 1 & \cdots & 1 \\
\alpha_{j_1} & \alpha_{j_2} & \cdots & \alpha_{j_m} \\
\alpha_{j_1}^2 & \alpha_{j_2}^2 & \cdots & \alpha_{j_m}^2 \\
\vdots & \vdots & \vdots & \vdots \\
\alpha_{j_1}^{m-1} & \alpha_{j_2}^{m-1} & \cdots & \alpha_{j_m}^{m-1} \\
\end{array}
\right)  
\left(
\begin{array}{c}
a_{j_1}e_{j_1} \\
a_{j_2}e_{j_2} \\
a_{j_3}e_{j_3} \\
\vdots \\
a_{j_m}e_{j_m} \\
\end{array}
\right). 
\]
The matrix $( \alpha_{j_t}^{s-1} )_{1 \le s, t \le m}$ is invertible, and hence 
the vector $a_{j_s} e_{j_s}$ can be written as a linear combination of 
$x, f(x), f^2(x), \ldots, f^{m-1}(x)$ for each $1 \le s \le m$. 
In particular, $e_i \in W$. This implies that $W \supseteq 
\oplus_{i \in I} k\cdot e_i$. So we have proved the lemma.    
\qed 

\begin{lemma}\label{lemma:existmatrix}
Let $V$ be a vector space over an infinite field $k$.  
For any non-zero vector  $v \in V$ and 
a finite subset $S \subset k^{\times}$, there exists $f \in {\rm GL}(V)$ satisfying  
the following conditions: 
\begin{enumerate}
\item There exists a basis $\{ v_1, v_2, \ldots, v_n \}$ of $V$ such that 
$v_i$ is an eigenvector of $f$ with eigenvalues $\beta_i \in k^{\times}\setminus S$ 
for $1 \le i \le n$. 
\item $\beta_1, \beta_2, \ldots, \beta_n$ are distinct. 
\item $v = v_1+v_2+\cdots +v_n$. 
\end{enumerate}  
In particular, $v$ is not contained in any proper $f$-invariant subspaces.  
\end{lemma}

{\it Proof}. 
Let us take vectors $v_1, v_2, \ldots, v_{n-1} \in V$ such that 
$\{ v, v_1, v_2, \ldots, v_{n-1} \}$ is a basis of $V$. 
Put $v_n := v - v_1 - v_2 - \cdots  - v_{n-1}$. 
Then $\{ v_1, v_2, \ldots, v_n \}$ is a basis of $V$ and  
$v = v_1+v_2+\cdots +v_n$. 
Let us choose distinct elements $\beta_1, \beta_2, \ldots, \beta_{n} \in k^{\times}\setminus S$. 
We define $f \in {\rm GL}(V)$ 
by $f(v_i) = \beta_i v_i$ for $1 \le i \le n$. By Lemma~\ref{lemma:linearalg}, for any proper 
$f$-invariant subspace $W$, there exists a proper subset $I$ of 
$\{ 1, 2, \ldots, n \}$ such that $W = \oplus_{i \in I} k\cdot v_i$.  
Hence $v = v_1+v_2+\cdots +v_n$ is not contained in $W$. 
This completes the proof. 
\qed

\begin{lemma}\label{lemma:eigenvector}
Let $k$ be a field. 
Let $A_1, A_2, \ldots, A_{\ell} \in {\rm GL}_m(k)$. 
Set $C := A_{\ell} A_{\ell -1} \cdots A_{2} A_{1}$ and \[
X = 
\left(
\begin{array}{cccccc}
0_m & 0_m & 0_m & \cdots & 0_m & A_{\ell} \\
A_1 & 0_m & 0_m & \cdots & 0_m & 0_m \\
0_m & A_2 & 0_m & \cdots & 0_m & 0_m \\
\vdots & \ddots & \ddots & \ddots & \vdots & \vdots \\
\vdots & \vdots & \ddots & \ddots & \ddots & \vdots \\
0_m & 0_m & 0_m & \cdots & A_{\ell-1} & 0_m \\
\end{array}
\right) 
\in {\rm GL}_n(k),  
\]
where $n=\ell m$ with $\ell \ge 2$. Suppose that the eigenvalues 
$\alpha_1, \ldots, \alpha_m$ of $C$ are distinct and that   
$\sharp \{ z \in k \mid z^{\ell} = \alpha_i \} = \ell$ for each  
$1 \le i \le m$. 
Then for each $\ell$-th root $\xi_{i, j}$ of $\alpha_i$  
$(1 \le i \le m, 1 \le j \le \ell)$ and for each eigenvector 
$v_i$ of $C$ with respect to $\alpha_i$, the vector 
%\[\begin{array}{l} 
\begin{multline*} 
w_{i, j} := \\  
{}^t ( \xi_{i, j}^{\ell -1} v_i,\:  \xi_{i, j}^{\ell -2} A_1 v_i,\: 
\xi_{i, j}^{\ell -3} (A_2 A_1) v_i, \ldots,\: 
\xi_{i, j} (A_{\ell - 2}\cdots A_{2} A_{1})v_i, (A_{\ell-1} A_{\ell-2} \cdots 
A_{2}A_{1})v_i )
\end{multline*} 
%\end{array}
%\]
is an eigenvector of $X$ with respect to the eigenvalue $\xi_{i, j}$. Conversely, all eigenvectors of $X$ can be 
obtained in this way (up to scalar multiplication). 
\end{lemma}

{\it Proof}. 
It is easy to check that $X w_{i, j} = \xi_{i, j}w_{i, j}$. 
The statement follows from that  $\{ \xi_{i, j} \mid 1 \le i \le m, 1 \le j \le \ell \}$ 
forms the set of $n$ distinct eigenvalues of $X$.  
\qed 

\bigskip

Let ${\rm F}_1 = \langle \alpha \rangle$ be the free group of rank $1$. 
By Proposition \ref{prop:3-2}, $X(d, n; {\rm F}_1)$ is a closed subscheme 
of ${\rm Rep}_n({\rm F}_1) \times {\rm Gr}(d, {\mathbb A}^n_{\mathbb Z})$. 
 Here recall that $X(d, n; {\rm F}_1) = \{ (\rho, W) \mid W \mbox{ is a } 
d\mbox{-dimensional } \rho(G)\mbox{-invariant 
subbundle of } {\Bbb A}^n \}$. 
Let $U(d, n) := U(d, n; {\rm F}_1)$ be the complement of 
$X(d, n; {\rm F}_1)$ in ${\rm Rep}_n({\rm F}_1) \times {\rm Gr}(d, {\mathbb A}^n_{\mathbb Z})$. Note that ${\rm Rep}_n({\rm F}_1) = {\rm GL}_n$ and that 
$U(d, n) = \{ (A, W) \mid W \mbox{ is not } A\mbox{-invariant} \} 
\subseteq {\rm GL}_n \times {\rm Gr}(d, {\mathbb A}^n_{\mathbb Z})$.
For a $X$-valued point $\phi$ of ${\rm Gr}(d, {\mathbb A}^n_{\mathbb Z})$ 
with a scheme $X$, 
denote by $\phi^{\ast}(W) \subset {\mathcal O}_X^{\oplus n}$  
the subbundle of rank $d$ induced by $\phi$ on $X$. 
Then we have the following diagram 
\[
\begin{array}{rcc}
{\rm GL}_{n, \phi} := \{ (A, x) \in {\rm GL}_n \times X \mid 
\phi^{\ast}(W)_{x} \mbox{ is not } A\mbox{-invariant} \} & 
\to & U(d, n) \\
\downarrow \;\;\;  & & \downarrow \\
{\rm GL}_n \times X & \stackrel{id \times \phi}{\to} & 
{\rm GL}_n \times {\rm Gr}(d, {\mathbb A}^n_{\mathbb Z}),  
\end{array}
\]
which is a fibre product. 
Hence ${\rm GL}_{n, \phi}$ is an open subscheme of 
${\rm GL}_n \times X$. 
 
In particular, for a geometric point $W$ of ${\rm Gr}(d, {\mathbb A}^{n}_{\Bbb Z})$, 
we have: 

\begin{proposition}\label{prop:GLW}
Let $k$ be an algebraically closed field. 
Let $W$ be a $k$-rational point of ${\rm Gr}(d, {\mathbb A}^{n}_{\Bbb Z})$. 
Then the subset $\{ A \in {\rm GL}_n(k) \mid W 
\mbox{ is not } A\mbox{-invariant } \}$ 
is an open subscheme of ${\rm GL}_n(k)$. 
\end{proposition} 

\begin{prop}\label{prop:existence}
Let $\ell$ and $m$ be positive integers with $\ell, m \ge 2$.  Set $n = \ell m$. 
Let $k$ be an algebraically closed field such that ${\rm ch}\: k$ does not divide $\ell$. 
Then there exists an irreducible representation $\rho : {\rm F}_2 \to {\rm GL}_n(k)$ 
of the free group ${\rm F}_2$ of rank $2$ such that $\wedge^m \rho$ has a realizable invariant subspace of dimension $\ell$. 
Moreover, there exists an irreducible representation $\rho : {\rm F}_2 \to {\rm GL}_n(k)$ 
such that $\rho$ is neither $m$-thick 
nor $\ell$-thick. 
\end{prop}

{\it Proof}. 
Let ${\rm F}_2 = \langle \alpha, \beta \rangle$. 
For constructing $\rho$, we need to determine 
$A := \rho(\alpha), B:=\rho(\beta) \in {\rm GL}_n(k)$. 
The group ${\rm GL}_n(k)$ acts canonically on $k^n$. 
Let $e_1, e_2, \ldots, e_n$ be the canonical basis 
of $k^n$.  
Set 
\[ 
A = 
\left(
\begin{array}{cccccc}
0_m & 0_m & 0_m & \cdots & 0_m & A' \\
I_m & 0_m & 0_m & \cdots & 0_m & 0_m \\
0_m & I_m & 0_m & \cdots & 0_m & 0_m \\
\vdots & \ddots & \ddots & \ddots & \vdots & \vdots \\
\vdots & \vdots & \ddots & \ddots & \ddots & \vdots \\
0_m & 0_m & 0_m & \cdots & I_m & 0_m \\
\end{array}
\right), 
B = 
\left(
\begin{array}{cccccc}
0_m & 0_m & 0_m & \cdots & 0_m & B_{\ell} \\
B_1 & 0_m & 0_m & \cdots & 0_m & 0_m \\
0_m & B_2 & 0_m & \cdots & 0_m & 0_m \\
\vdots & \ddots & \ddots & \ddots & \vdots & \vdots \\
\vdots & \vdots & \ddots & \ddots & \ddots & \vdots \\
0_m & 0_m & 0_m & \cdots & B_{\ell-1} & 0_m \\
\end{array}
\right),   
\]
where $A', B_1, \ldots, B_{\ell} \in {\rm GL}_{m}(k)$ 
will be suitably chosen. Let us define 
${\Phi} : {\rm GL}_m(k) \times \cdots \times {\rm GL}_m(k) = {\rm GL}_m(k)^{\ell} 
\to {\rm GL}_n(k)$ by $(B_1, B_2, \ldots, B_{\ell}) \mapsto B$.

First, we show that $\rho$ is not $m$-thick.
Let $W := \langle e_1 \wedge e_2 \wedge \cdots \wedge e_m, 
e_{m+1}\wedge \cdots \wedge e_{2m}, 
e_{2m+1} \wedge \cdots \wedge e_{3m}, 
\ldots, 
e_{(\ell-1)m+1} \wedge \cdots \wedge e_{n}
\rangle 
\subseteq \wedge^{m} V$.  
Note that $\wedge^m \rho$ has a realizable invariant subspace $W$ of 
dimension $\ell$.  
Since $\ell \le m(n-m)$, $\rho$ is not $m$-thick 
by Proposition \ref{prop:condnothick}.  

Second, we show that $\rho$ is not $\ell$-thick. 
For $1 \le i \le \ell$, we put $J_{i} := \{ (i-1)m+1, 
(i-1)m+2, \ldots, im \}$. 
Then $J_{1} \sqcup J_{2} \sqcup \cdots \sqcup J_{\ell} = 
\{ 1, 2, \ldots, n \}$. 
Let $Y$ be the subspace of 
$\wedge^{\ell} V$ generated by 
$\{ e_{i_1} \wedge e_{i_2} \wedge \cdots \wedge e_{i_{\ell}} 
\mid i_1 \in J_1, i_2 \in J_2, \ldots, i_{\ell} \in J_{\ell}  \}$. 
Note that $Y$ is an $m^{\ell}$-dimensional $(\wedge^{\ell} \rho) ({\rm F}_2)$-invariant  realizable subspace of $\wedge^{\ell} V$. 
The subspace $Y^{\perp}$ of 
$\wedge^{n-\ell} V$ contains  
$e_{1} \wedge e_{2} \wedge \cdots \wedge e_{m} \wedge v'$ 
for any $v' \in \wedge^{n - \ell - m} V$. 
In particular, $Y^{\perp}$ is realizable. 
By Proposition~\ref{prop:condofm-thick}, $\rho$ is not $\ell$-thick.

Finally, we show that $\rho$ is irreducible if $A', B_1, \ldots, B_{\ell}$ are 
suitably chosen. 
Let $A' = {\rm diag} ( \alpha_1, \alpha_2, \ldots, \alpha_m )$, where 
$\alpha_1, \ldots, \alpha_m \in k^{\times}$ are distinct.  
For each $\ell$-th root $\xi_{i, j}$ of $\alpha_i$ $(1 \le i \le m, 1 \le j \le \ell)$, 
we define $w_{i, j} := {}^t (\xi_{i, j}^{\ell-1} e'_i, \xi_{i, j}^{\ell-2} e'_i, 
\ldots, e'_i)$ as in Lemma \ref{lemma:eigenvector}, where 
$A_1 = A_2 = \cdots A_{\ell-1} = I_m$ and $A_{\ell}=A'$. 
Here we use $e'_1, \ldots, e'_m$ as the canonical basis of $k^m$ in the sequel. 
Then $w_{i, j}$ 
is an eigenvector of $A$.  

By Lemma \ref{lemma:linearalg}, for any $A$-invariant subspace $W$ of $k^n$, there 
exists a subset $I$ of $\{ (i, j) \mid 1 \le i \le m, 1 \le j \le \ell \}$
such that $W = W_{I} := \oplus_{(i, j) \in I} k\cdot w_{i, j}$. 
For proving $\rho$ is irreducible, it suffices to show that 
$B$ does not keep any non-trivial $A$-invariant subspace $W_{I}$ invariant.  
For each non-trivial $A$-invariant subspace $W_{I}$, 
we set ${\rm GL}_n(k)_{I} := \{ B \in {\rm GL}_n(k) 
\mid W_I \mbox{ is not $B$-invariant } \}$. 
By Proposition~\ref{prop:GLW}, ${\rm GL}_n(k)_{I}$ is an open subscheme 
of ${\rm GL}_n(k)$. Let us prove the claim that the open subset  
$\Phi^{-1}({\rm GL}_n(k)_{I}) \subseteq {\rm GL}_m(k)^\ell $ 
is not empty for each non-empty 
proper subset $I$ of $\{ (i, j) \mid 1 \le i \le m, 1 \le j \le \ell \}$. 
If $\Phi^{-1}({\rm GL}_n(k)_{I}) \neq \emptyset$ for each $I$, then 
$\cap_{I} \Phi^{-1}({\rm GL}_n(k)_{I}) \neq \emptyset$ because  
${\rm GL}_m(k)^{\ell}$ is irreducible. 
Then by taking $(B_1, \ldots, B_{\ell}) \in { \cap_{I}} 
\Phi^{-1}({\rm GL}_n(k)_{I})$, 
we obtain an irreducible representation $\rho$, which completes the proof. 

For proving the claim that $\Phi^{-1}({\rm GL}_n(k)_{I}) \neq \emptyset$, 
take some $(i_0, j_0) \in I$.   
By Lemma~\ref{lemma:existmatrix}, there exist a basis $\{v_1, v_2, \ldots, v_m \}$ 
of $k^m$ and $f : k^m \to k^m$ such that $f(v_i) = \beta_i v_i$ ($1 \le i \le m$) and 
$e'_{i_0} = v_1+\cdots +v_m$.  
Here $\beta_1, \ldots, \beta_m$ are distinct elements in $k^{\times} \setminus \{  
\alpha_{i}  \mid 1 \le i \le m \}$.  Let $B_{\ell} \in {\rm GL}_m(k)$ 
be the corresponding matrix to $f$. 
Set $B_1 = B_2 = \cdots = B_{\ell -1} = I_m$ and $B = \Phi(B_1, \ldots, B_{\ell})$. 
For each $\ell$-th root $\eta_{ij}$ of $\beta_i$ ($1 \le i \le m, 1 \le j \le \ell$),  
put $w'_{i, j} := {}^t (\eta^{\ell-1}_{i, j} v_i, \eta^{\ell-2}_{i, j} v_i, 
\ldots, v_i) \in k^n$.  By Lemma \ref{lemma:eigenvector},   
$B w'_{i,j} = \eta_{i, j} w'_{i, j}$ for each $i, j$.   
Since $\{ w'_{i, j} \mid 1 \le i \le m, 1 \le j \le \ell \}$ is a basis of $k^n$, 
we can write $w_{i_0, j_0} = \sum c_{i, j} w'_{i, j}$ for $c_{i, j} \in k$.  
If $c_{i, j} \neq 0$ for all $i, j$, then $w_{i_0, j_0} \in W_{I}$ is not 
contained in any non-trivial $B$-invariant subspaces 
by Lemma~\ref{lemma:linearalg}.    
In particular, $W_{I}$ is not $B$-invariant and 
$(B_1, \ldots, B_{\ell}) \in \Phi^{-1}({\rm GL}_n(k)_{I})$, which implies the claim.  
Hence we only need to show that $c_{i, j} \neq 0$ for all $i, j$. 

Let us show that $c_{i, j} \neq 0$. 
For each $1 \le i \le m$, we define the $\ell$-dimensional subspace $U_i := 
\langle {}^{t} (v_i, 0, 0, \ldots, 0), {}^{t} (0, v_i, 0, \ldots, 0), 
\ldots, {}^{t} (0, 0, 0, \ldots, v_{i}) \rangle \subset k^n$. 
Let $p_i : k^n = U_1 \oplus \cdots \oplus U_m \to U_i$ be the 
projection onto $U_i$. 
Since $U_i = \oplus_{1 \le j \le \ell} \:k\cdot w'_{i, j}$,  
\[
p_i(w_{i_0, j_0}) = p_i({}^t (\xi_{i_0, j_0}^{\ell-1} e'_{i_0}, \xi_{i_0, j_0}^{\ell-2} e'_{i_0}, 
\ldots, e'_{i_0}))=\sum_{1 \le j \le \ell} c_{i, j} w'_{i, j}.
\]  
On the other hand, $e'_{i_0} =  v_1+\cdots +v_m$ and hence 
\[
{}^t (\xi_{i_0, j_0}^{\ell-1} v_{i}, \xi_{i_0, j_0}^{\ell-2} v_{i}, 
\ldots, v_{i}) =\sum_{1 \le j \le \ell} c_{i, j} w'_{i, j}.
\]  
Then we have 
\[
\left( 
\begin{array}{cccc} 
\eta^{\ell-1}_{i,1} & \eta^{\ell-1}_{i,2} & \cdots & \eta^{\ell-1}_{i,\ell} \\
\eta^{\ell-2}_{i,1} & \eta^{\ell-2}_{i,2} & \cdots & \eta^{\ell-2}_{i,\ell} \\
\vdots & \vdots & \cdots & \vdots \\
1 & 1 & \cdots & 1 \\ 
\end{array} 
\right) 
\left( 
\begin{array}{c} 
c_{i,1} \\
c_{i,2} \\
\vdots \\
c_{i,\ell} 
\end{array}  
\right) 
= 
\left( 
\begin{array}{c} \vspace*{0.5ex} 
\xi^{\ell-1}_{i_0,j_0} \\ 
\xi^{\ell-2}_{i_0,j_0} \\
\vdots \\
1  
\end{array}  
\right). 
\]
By Cramer's rule, 
\[ 
c_{i,j} = \det \left( 
\begin{array}{cccccc} \vspace*{1ex}
\eta^{\ell-1}_{i,1} & \eta^{\ell-1}_{i,2} & \cdots & \overset{j}{\check{\xi_{i_0, j_0}^{\ell-1}}} 
& \cdots & \eta^{\ell-1}_{i,\ell} \\
\eta^{\ell-2}_{i,1} & \eta^{\ell-2}_{i,2} & \cdots & \xi_{i_0, j_0}^{\ell-2} & 
\cdots & \eta^{\ell-2}_{i,\ell} \\
\vdots & \vdots & \cdots & \vdots & \cdots & \vdots \\
1 & 1 & \cdots & 1 & \cdots & 1 \\ 
\end{array} 
\right) \cdot 
\det \left( 
\begin{array}{cccc} 
\eta^{\ell-1}_{i,1} & \eta^{\ell-1}_{i,2} & \cdots & \eta^{\ell-1}_{i,\ell} \\
\eta^{\ell-2}_{i,1} & \eta^{\ell-2}_{i,2} & \cdots & \eta^{\ell-2}_{i,\ell} \\
\vdots & \vdots & \cdots & \vdots \\
1 & 1 & \cdots & 1 \\ 
\end{array} 
\right)^{-1}.  
\]
The Vandermonde determinant is not $0$ because 
$\eta_{i,j}$ and $\xi_{i_0, j_0}$ are distinct. Hence $c_{i,j} \neq 0$.   
Therefore we have completed the proof.
\hfill $\Box$ 

\begin{corollary}\label{cor:thdivide}
If $m$ divides $n$, then $r({\wedge^m (n)}) = [\frac{n-1}{m}] + 1 = \frac{n}{m}$. 
\end{corollary}

{\it Proof}. 
If $m=1$ or $n=m$, then the statement is trivial. 
Let $n > m \ge 2$. 
By Proposition \ref{prop:existence},   
there exists an $n$-dimensional irreducible representation $\rho$ of ${\rm F}_2$ 
such that $\wedge^{m} \rho$ has a realizable 
invariant subspace of dimension $n/m$. 
Hence we have 
$r({\wedge^m (n)}) = n/m$ by Proposition \ref{prop:rinequality}. 
\hfill $\Box$

\bigskip 

The following remark was suggested by the referee. 

\begin{remark}\label{remark:f3rep}\rm
Proposition \ref{prop:existence} shows that there exists an irreducible representation $\rho : {\rm F}_2 \to {\rm GL}_n(k)$ such that $\rho$ is neither $m$-thick nor $\ell$-thick, where $n=\ell m$ for $\ell, m \ge 2$. 
However, we can easily construct an 
irreducible representation $\rho : {\rm F}_3 \to {\rm GL}_n(k)$ 
of the free group ${\rm F}_3 = \langle \alpha, \beta, \gamma \rangle$ of rank $3$   
such that $\rho$ is neither $m$-thick nor $\ell$-thick. Indeed, set 
$\rho(\alpha) = {\rm diag}(\alpha_1, \alpha_2, \ldots, \alpha_n)$,  
\[
\rho(\beta) = \left(
\begin{array}{cc}
B' & 0_{m, n-m} \\
0_{n-m, m} & I_{n-m} \\ 
\end{array}
\right), \; 
\rho(\gamma) = 
\left(
\begin{array}{cccccc}
0_m & 0_m & 0_m & \cdots & 0_m & I_m \\
I_m & 0_m & 0_m & \cdots & 0_m & 0_m \\
0_m & I_m & 0_m & \cdots & 0_m & 0_m \\
\vdots & \ddots & \ddots & \ddots & \vdots & \vdots \\
\vdots & \vdots & \ddots & \ddots & \ddots & \vdots \\
0_m & 0_m & 0_m & \cdots & I_m & 0_m \\
\end{array}
\right), 
\] 
where $\alpha_1, \ldots, \alpha_n \in k^{\times}$ are distinct and $B' = \left(
\begin{array}{cccccc}
0 & 0 & 0 & \cdots & 0 & 1 \\
1 & 0 & 0 & \cdots & 0 & 0 \\
0 & 1 & 0 & \cdots & 0 & 0 \\
\vdots & \ddots & \ddots & \ddots & \vdots & \vdots \\
\vdots & \vdots & \ddots & \ddots & \ddots & \vdots \\
0 & 0 & 0 & \cdots & 1 & 0 \\
\end{array}
\right) \in {\rm GL}_m(k)$. 
We can show that $\rho : {\rm F}_3 \to {\rm GL}_n(k)$ is irreducible by Lemma \ref{lemma:linearalg}.   
We can also prove that $\rho$ is neither $m$-thick nor $\ell$-thick, and that 
$W := \langle e_1 \wedge e_2 \wedge \cdots \wedge e_m, 
e_{m+1}\wedge \cdots \wedge e_{2m}, 
e_{2m+1} \wedge \cdots \wedge e_{3m}, 
\ldots, 
e_{(\ell-1)m+1} \wedge \cdots \wedge e_{n}
\rangle 
\subseteq \wedge^{m} V$ is a realizable $(\wedge^m \rho)({\rm F}_3)$-invariant subspace of dimension $\ell$ 
as in the proof of Proposition \ref{prop:existence}. 
\end{remark}

\bigskip 

By the definition, it is obvious that  
$r(\wedge^m(n)) \le \binom{n}{m}$. 
The following proposition gives us a non-trivial upper bound of $r(\wedge^m(n))$. 

\begin{prop}\label{prop:r-upper} 
For $0 < m <n$, $r(\wedge^m(n)) \le n$. 
\end{prop} 

{\it Proof}. 
Let $a$ and $b$ are distinct non-zero elements of  a field $k$. Assume that 
$\sharp \{ c\in k \mid c^n = a \} = \sharp \{ c\in k \mid c^n = b \} = n$. 
Let $\{ \xi_{i} \mid 1 \le i \le n \}$ and $\{ \eta_{i} \mid 1 \le i \le n\}$ be 
the $n$-th roots of $a$ and $b$, respectively. 
We define an $n$-dimensional representation $\rho$ of 
the free group ${\rm F}_2 = \langle \alpha, \beta \rangle$ by 
\[ 
\rho(\alpha) = 
\left(
\begin{array}{cccccc}
0 & 0 & 0 & \cdots & 0 & a \\
1 & 0 & 0 & \cdots & 0 & 0 \\
0 & 1 & 0 & \cdots & 0 & 0 \\
\vdots & \ddots & \ddots & \ddots & \vdots & \vdots \\
\vdots & \vdots & \ddots & \ddots & \ddots & \vdots \\
0 & 0 & 0 & \cdots & 1 & 0 \\
\end{array}
\right), \;\; 
\rho(\beta) = 
\left(
\begin{array}{cccccc}
0 & 0 & 0 & \cdots & 0 & b \\
1 & 0 & 0 & \cdots & 0 & 0 \\
0 & 1 & 0 & \cdots & 0 & 0 \\
\vdots & \ddots & \ddots & \ddots & \vdots & \vdots \\
\vdots & \vdots & \ddots & \ddots & \ddots & \vdots \\
0 & 0 & 0 & \cdots & 1 & 0 \\
\end{array}
\right).    
\]
Set $v_i = {}^t(\xi_{i}^{n-1}, \xi_{i}^{n-2}, \ldots, \xi_{i}, 1)$ and $w_i = {}^t(\eta_{i}^{n-1}, \eta_{i}^{n-2}, \ldots, \eta_{i}, 1)$. 
It is easy to check that $\rho(\alpha)v_i = \xi_i v_i$ and $\rho(\beta)w_i=\eta_i w_i$  
for $1 \le i \le n$ and that $\{ v_1, \ldots, v_n \}$ and $\{ w_1, \ldots, w_n \}$ are bases of $k^n$.   
In a similar way as the last part of the proof of Proposition~\ref{prop:existence}, 
we can prove that $\rho : {\rm F}_2 \to {\rm GL}_n(k)$ is irreducible. 
Indeed, let $V$ be a non-zero subspace of $k^n$ which is invariant under $\rho(\alpha)$ and $\rho(\beta)$. 
By Lemma \ref{lemma:linearalg}, there exist subsets $I$ and $J$ of  
$\{ 1, 2, \ldots, n \}$ such that $V= \oplus_{i \in I} k \cdot v_i = \oplus_{j \in J} k \cdot w_j$.    
Suppose that $w_j \in V$.  Put $w_j = \sum c_i v_i$. Then we have 
\[
\left( 
\begin{array}{cccc} 
\xi^{n-1}_{1} & \xi^{n-1}_{2} & \cdots & \xi^{n-1}_{n} \\
\xi^{n-2}_{1} & \xi^{n-2}_{2} & \cdots & \xi^{n-2}_{n} \\
\vdots & \vdots & \cdots & \vdots \\
1 & 1 & \cdots & 1 \\ 
\end{array} 
\right) 
\left( 
\begin{array}{c} 
c_{1} \\
c_{2} \\
\vdots \\
c_{n} 
\end{array}  
\right) 
= 
\left( 
\begin{array}{c} \vspace*{0.5ex} 
\eta^{n-1}_{j} \\ 
\eta^{n-2}_{j} \\
\vdots \\
1  
\end{array}  
\right). 
\]
By Cramer's rule, 
\[ 
c_{i} = \det \left( 
\begin{array}{cccccc} \vspace*{1ex}
\xi^{n-1}_{1} & \xi^{n-1}_{2} & \cdots & \overset{i}{\check{\eta_{j}^{n-1}}} 
& \cdots & \xi^{n-1}_{n} \\
\xi^{n-2}_{1} & \xi^{n-2}_{2} & \cdots & \eta_{j}^{n-2} & 
\cdots & \xi^{n-2}_{n} \\
\vdots & \vdots & \cdots & \vdots & \cdots & \vdots \\
1 & 1 & \cdots & 1 & \cdots & 1 \\ 
\end{array} 
\right) \cdot 
\det \left( 
\begin{array}{cccc} 
\xi^{n-1}_{1} & \xi^{n-1}_{2} & \cdots & \xi^{n-1}_{n} \\
\xi^{n-2}_{1} & \xi^{n-2}_{2} & \cdots & \xi^{n-2}_{n} \\
\vdots & \vdots & \cdots & \vdots \\
1 & 1 & \cdots & 1 \\ 
\end{array} 
\right)^{-1}.  
\]
The Vandermonde determinant is not $0$ because 
$\xi_{i}$ and $\eta_{j}$ are distinct, and hence $c_i \neq 0$ for $1 \le i \le n$. We see that  $I = \{ 1, 2, \ldots, n \}$ and that $V = k^n$.   This implies that $\rho$ is irreducible.

Let $e_1, e_2, \ldots, e_n$ be the canonical basis of $k^n$. 
For $0 < m < n$, define an $n$-dimensional subspace $W_m$ of $\wedge^m k^n$ by  
\begin{multline*} 
W_m := \langle e_1\wedge e_2 \wedge \cdots \wedge e_m,  
e_2 \wedge e_3 \wedge \cdots \wedge e_{m+1}, 
\ldots, e_{i}\wedge e_{i+1} \wedge \cdots \wedge e_{i+m-1}, \\ \ldots, 
e_{n-1} \wedge e_{n} \wedge e_1 \wedge \cdots \wedge e_{m-2}, 
e_{n} \wedge e_1 \wedge \cdots \wedge e_{m-1} 
\rangle \subset \wedge^m k^n.
\end{multline*} 
Then $W_m$ is an ${\rm F}_2$-invariant realizable subspace of $\wedge^m k^n$. 
Hence $r(\wedge^m(n)) \le \dim W_m = n$. This completes the proof. 
\hfill $\Box$

\bigskip 

We prepare some basic results on perfect pairings for 
determining some $r(\wedge^m(n))$.  
In the sequel, by a $G$-module we understand a finite-dimensional   
left $G$-module over a field $k$ 
for a group $G$. For a $G$-module $W$, the dual $W^{\ast}$ 
is defined as $W^{\ast} := \{ f : W \to k \mid k \mbox{-linear } \}$, where 
$(g\cdot f)(\ast) := f(g^{-1} \ast)$ for $g \in G$ and $f \in W^{\ast}$.   
For $G$-modules $W, W'$, we define the $G$-module $W \otimes_{k} W'$
by $g\cdot (u\otimes v) := gu \otimes gv$ for $g \in G$, 
$u \in W$, and $v \in W'$.
 
\begin{lemma}\label{lemma:perfect}
Let $W, W'$ be finite-dimensional $G$-modules over a field $k$. 
Let $L$ be a one-dimensional $G$-module. 
Suppose that $f : W \times W' \to L$ is a $G$-equivariant 
perfect pairing. In other words, the bilinear map $f$ satisfies 
\begin{enumerate}
\item $f(u, v)=0$ for all $v \in W'  \Rightarrow u=0$, 
\item $f(u, v)=0$ for all $u \in W   \Rightarrow v=0$, 
\item $f(gu, gv)=g(f(u, v))$ for all $g \in G, u \in W, v \in W'$.   
\end{enumerate}  
Then there exists a canonical isomorphism $W' \cong W^{\ast}\otimes_{k} L$ 
as $G$-modules. 
\end{lemma}

{\it Proof}.
Let $e$ be a non-zero vector of $L$.  
Let $\phi_{e} : k \to L$ be 
the linear isomorphism defined by $a \mapsto ae$ for $a \in k$. 
We define the linear map $\Phi : W' \to W^{\ast} \otimes L$ 
by $v \mapsto \phi_{e}^{-1}(f(\ast, v)) \otimes e$. 
Note that the definition of $\Phi$ is independent from the choice of $e$. 
We claim that $\Phi$ is an isomorphism as $G$-modules.

First, we show that $\Phi$ is $G$-equivariant. 
Let $\chi : G \to {\rm GL}_1(k)$ be the character associated to $L$.  
In other words, $g\cdot w = \chi(g) w$ for $g \in G$ and $w \in L$. 
We see that 
\begin{eqnarray*}
\Phi(gv)  & = & \phi_{e}^{-1}(f(\ast, gv))\otimes e = \phi_{e}^{-1}(g\cdot (f(g^{-1}\ast, v)))\otimes e
= \phi_{e}^{-1}(\chi(g) f(g^{-1} \ast, v)) \otimes e \\
 & = & \phi_{e}^{-1}( f(g^{-1}\ast, v))\otimes \chi(g) e 
= \phi_{e}^{-1}( f(g^{-1}\ast, v))\otimes g\cdot e = g\cdot \Phi(v).  
\end{eqnarray*} 
Hence $\Phi$ is $G$-equivariant. 

Next, suppose that $\Phi(v)=0$. The assumption implies that 
$f(u, v)=0$ for all $u \in W$. Because of perfectness, we have $v=0$. 
Thus we proved that $\Phi$ is injective. 
On the other hand, we see that $\dim W' = \dim (W^{\ast}\otimes L)$, 
which implies that $\Phi$ is surjective. 
Therefore $\Phi$ is an isomorphism.  
\hfill $\Box$ 

\begin{corollary}\label{cor:perfect} 
Let $W, W'$ be finite-dimensional $G$-modules over a field $k$. 
Let $L$ be a one-dimensional $G$-module. 
Suppose that a 
bilinear map $f : W \times W' \to L$ satisfies:  
\begin{enumerate}
\item $f(u, v)=0$ for all $v \in W'  \Rightarrow u=0$.  
\item $f(gu, gv)=g(f(u, v))$ for all $g \in G, u \in W, v \in W'$.   
\end{enumerate}  
Then there exists a canonical surjection $W' \to W^{\ast}\otimes_{k} L$ as $G$-modules. 
\end{corollary}

{\it Proof}.
Let $W'^{\sharp} := \{ v \in W' \mid f(u, v)=0 \mbox{ for all }u \in W \}$. 
The bilinear map $f : W \times W' \to L$ induces a $G$-equivariant perfect 
pairing $\overline{f} : W \times (W'/W'^{\sharp}) \to L$. 
By Lemma \ref{lemma:perfect} we have a canonical isomorphism 
$\Phi : (W'/W'^{\sharp}) \cong W^{\ast}\otimes L$. 
Composing $\Phi$ and the projection $W' \to (W'/W'^{\sharp})$, we have a 
canonical surjection $W' \to W^{\ast}\otimes_{k} L$. \\  
\hfill $\Box$ 

\begin{corollary}\label{cor:perfperp} 
Let $W, W'$ be finite-dimensional $G$-modules over a field $k$. 
Let $L$ be a one-dimensional $G$-module. 
Let $Z$ be an irreducible $G$-submodule of $W$, and let 
$Y$ be a $G$-submodule of $W'$. 
Suppose that any $G$-homomorphism 
$\phi : Y \to Z^{\ast}\otimes L$ is not surjective. 
If $f : W \times W' \to L$ is a $G$-equivariant 
perfect pairing, then  
$f(z, y)=0$ for all $z \in Z, y \in Y$. 
\end{corollary}

{\it Proof}.
Let $Y^{\sharp} := \{ y \in Y \mid f(z, y)=0 \mbox{ for all }z \in Z \}$. 
If $Y^{\sharp} = Y$, then the statement is true. 
Suppose that $Y^{\sharp} \neq Y$. 
Then $f$ induces $\overline{f} : Z \times (Y/Y^{\sharp}) \to L$ 
which has the property that $\overline{f}(z, \overline{y})=0$ 
for all $z \in Z$ implies $\overline{y} = 0$. 
By Corollary \ref{cor:perfect}, 
there exists a surjection $\phi : Z \to (Y/Y^{\sharp})^{\ast}\otimes L$. 
Since $Z$ is irreducible, $\phi$ is an isomorphism. 
Taking $\phi \otimes L^{\ast}$ and the dual, we have $Z^{\ast} \otimes L 
\cong (Y/Y^{\sharp})$.  
Then we obtain a surjection 
$Y \to (Y/Y^{\sharp}) \cong Z^{\ast}\otimes L$, which is a contradiction. 
Hence $Y^{\sharp} = Y$. \\ 
\hfill $\Box$ 

\begin{proposition}\label{prop:rdual}
For $0 < m < n$, $r(\wedge^m(n)) = r(\wedge^{n-m}(n))$. 
\end{proposition} 

{\it Proof}. 
Let $\rho : G \to {\rm GL}(V)$ be an $n$-dimensional irreducible 
representation of a group $G$. 
Assume that $\wedge^m V$ has a $G$-invariant realizable subspace $W$ 
of dim $d$. 
We claim that $\wedge^{n-m}(V^{\ast})$ has a $G$-invariant 
realizable subspace of dim $d$. 
Since $V^{\ast}$ is an $n$-dimensional irreducible $G$-module, 
we see that $r(\wedge^m(n)) \ge r(\wedge^{n-m}(n))$ from this claim. 
By Changing $m$ and $n-m$, we have $r(\wedge^{n-m}(n)) \ge r(\wedge^{m}(n))$, and we can conclude that 
$r(\wedge^m(n)) = r(\wedge^{n-m}(n))$. 

Let us prove the claim. 
Considering the perfect pairing $\wedge^mV \times 
\wedge^{n-m} V \to \wedge^nV$,  
we have  a canonical isomorphism $\Phi : 
\wedge^m V \cong (\wedge^{n-m} V)^{\ast} \otimes \wedge^{n} V$ 
by Lemma \ref{lemma:perfect}. 
Let $e_1, \ldots, e_n$ be a basis of $V$ such that 
$e_1 \wedge \cdots \wedge e_m \in W$. 
Let $f_1, \ldots, f_n$ be the dual basis for $e_1, \ldots, e_n$. 
Set $W' := \Phi(W) \otimes (\wedge^n V)^{\ast}$. Then 
$W'$ is a $d$-dimensional $G$-invariant subspace 
of $(\wedge^{n-m} V)^{\ast} \otimes \wedge^{n} V \otimes (\wedge^n V)^{\ast} 
\cong \wedge^{n-m}(V^{\ast})$. 
We easily see that $W'$ contains $(e_1 \wedge \cdots \wedge e_{m}) \wedge \ast 
= f_{m+1}\wedge \cdots \wedge f_n$. 
This implies that 
$\wedge^{n-m}(V^{\ast})$ has a $G$-invariant 
realizable subspace $W'$ of dim $d$. 
Therefore we have proved the statement.    
\hfill $\Box$ 

\begin{remark}\rm
By the definition, $r(\wedge^0(n)) = r(\wedge^n(n))=1$. 
Hence $r(\wedge^m(n)) = r(\wedge^{n-m}(n))$ for $0 \le m \le n$. 
%By Corollaries~\ref{cor:rge2} and \ref{cor:thdivide} and Proposition 
%\ref{prop:r-upper}, 
We see that 
$r(\wedge^{1}(n)) = r(\wedge^{n-1}(n))=n$ for $n \ge 2$.  
By Corollary \ref{cor:rge2} and Proposition \ref{prop:r-upper},  
$2 \le r(\wedge^{m}(n)) \le n$ for $0 < m < n$. 
It is not easy to calculate the $r$-number $r(\wedge^m(n))$ in general. 
\end{remark}

\begin{proposition}\label{prop:conj5}
$r(\wedge^2(5)) = r(\wedge^3(5)) \ge 4$. 
\end{proposition}

{\it Proof}. 
Since $r(\wedge^2(5)) = r(\wedge^3(5))$ by Proposition \ref{prop:rdual}, 
it suffices to prove that $r(\wedge^2(5)) \ge 4$. 
By Proposition \ref{prop:rinequality} we have $r(\wedge^2(5)) \ge 3$. 
We claim that $r(\wedge^2(5)) \neq 3$. 
Suppose that there exists a $3$-dimensional realizable invariant subspace $W$ of  
$\wedge^2 V$ for a $5$-dimensional irreducible representation 
$\rho : G \to {\rm GL}(V)$ of a group $G$. 
Since $W$ is realizable, there exist linearly independent vectors 
$e_1, e_2 \in V$ such that $e_1 \wedge e_2 \in W$. 
By irreducibility of $\rho$, 
there exists $g \in G$ such that $\rho(g)(e_1)$ can not be written 
as a linear combination of $\{ e_1, e_2 \}$. 
Similarly, there exists $g' \in G$ such that $\rho(g')(e_1)$ can not be written 
as a linear combination of $\{ e_1, e_2, \rho(g)e_1, \rho(g)e_2 \}$. 
Put $v_1 := e_1 \wedge e_2$, $v_2 := \rho(g)e_1 \wedge \rho(g)e_2$, and 
$v_3 := \rho(g')e_1 \wedge \rho(g')e_2$. 
Note that $\{ e_1, e_2, \rho(g)e_1, \rho(g')e_1 \}$ and 
$\{ \rho(g)e_1, \rho(g)e_2, \rho(g')e_1 \}$ are linearly independent. 
Then $v_1 \wedge \rho(g)e_1 \wedge \rho(g')e_1 \neq 0$ and 
$v_2 \wedge \rho(g')e_1 \neq 0$.  
We easily see that $v_1, v_2, v_3 \in W$ are linearly independent. 
Hence $W = \langle v_1, v_2, v_3 \rangle$. 

We define the subspace $W \wedge W$ of $\wedge^4 V$ as the subspace spanned by the vectors  
$\{ x \wedge y \in \wedge^4 V \mid x, y \in W \}$. 
The vector space $W \wedge W$ can be spanned by the vectors 
$v_1 \wedge v_2$, $v_1 \wedge v_3$, and $v_2 \wedge v_3$. 
Hence $W \wedge W$ is a $G$-invariant subspace of $\wedge^4 V$ and $\dim W \wedge W \le 3$. 
Since $V$ is irreducible, $\wedge^4 V$ is also irreducible. 
Thus $W \wedge W = 0$. 
Then there exist $g_1, g_2, g_3 \in G$ such that $\{ e_1, e_2, \rho(g_1)e_1, \rho(g_2)e_1, \rho(g_3)e_1 \}$ 
is linearly independent and $(e_1 \wedge e_2) \wedge (\rho(g_i)e_1 \wedge \rho(g_i)e_2)=0$ 
for $1 \le i \le 3$. 
The vector $\rho(g_i)e_2$ can be written as a linear combination of 
$\{ e_1, e_2, \rho(g_i)e_1 \}$ for each $i$. 
So we easily see that 
$e_1 \wedge e_2$, $\rho(g_1)e_1 \wedge \rho(g_1)e_2$,  
$\rho(g_2)e_1 \wedge \rho(g_2)e_2$, and $\rho(g_3)e_1 \wedge \rho(g_3)e_2$ 
are linearly independent. 
Thus $\dim W \ge 4$. This is a contradiction.
Therefore  
$r(\wedge^2(5)) \ge 4$. 
\hfill $\Box$

\bigskip 
Later we will show that  
$r(\wedge^2(5)) = r(\wedge^3(5)) = 4$ in Proposition~\ref{prop:r-number2-5}. 
%\bigskip 

\begin{proposition}
$r(\wedge^2(6))=r(\wedge^4(6))=3$ and $r(\wedge^3(6))=2$.  
\end{proposition}

{\it Proof}.
By Corollary~\ref{cor:thdivide} and 
Proposition~\ref{prop:rdual}, we can verify 
the statement. 
\hfill $\Box$ 

%%%%%%%%%%%%%%%%%%%%%%%%%%%%%%%%%%%%%%%%%%%%%%%%%%%%%%%%%%%%%%%%%%%%%%%%%%%%%%%%%%%%%%%55

\section{Criterion for thickness}

In this section, we discuss criteria for thickness. 
First, we deal with $4$-dimensional representations.

\begin{proposition}\label{prop:eq4}
Let $V$ be a $4$-dimensional vector space over 
an algebraically closed field $k$. 
For a representation $\rho : G \to {\rm GL}(V)$, 
the following statements are equivalent: 
\begin{enumerate}
\item\label{cond4-1} 
$\rho$ is thick.  
\item\label{cond4-2} 
$\rho$ is $2$-thick.  
\item\label{cond4-3} 
$\rho$ is irreducible and the induced representation 
$\wedge^2 \rho : G \to {\rm GL}(\wedge^2 V)$ has no $G$-invariant 
subspace $W \subset \wedge^2 V$ such that $2 \le \dim W \le 4$.  
\item\label{cond4-4} 
$\rho$ is irreducible and the induced representation 
$\wedge^2 \rho : G \to {\rm GL}(\wedge^2 V)$ has no $G$-invariant 
subspace $W \subset \wedge^2 V$ such that $\dim W = 2$ or $3$.  
\end{enumerate}
\end{proposition}

{\it Proof}. 
It is trivial that $(\ref{cond4-1}) \Rightarrow (\ref{cond4-2})$ and 
$(\ref{cond4-3}) \Rightarrow (\ref{cond4-4})$. 
If $\rho$ is $2$-thick, then $\rho$ is irreducible by Proposition 
\ref{prop:imp}. Hence by Proposition~\ref{prop:dual},  
$\rho$ is $m$-thick for $1 \le m \le 4$, which 
implies that $(\ref{cond4-2}) \Rightarrow (\ref{cond4-1})$. 
Assume that $\rho$ satisfies $(\ref{cond4-4})$. 
Suppose that $\wedge^2 \rho$ has a $G$-invariant subspace $W \subset \wedge^2 V$ 
of $\dim W=4$. Then $W^{\perp} \subset \wedge^2 V$ is a $2$-dimensional  
$G$-invariant subspace. This is a contradiction. Hence 
$(\ref{cond4-4}) \Rightarrow (\ref{cond4-3})$. 
   
Next, we show that $(\ref{cond4-2}) \Rightarrow (\ref{cond4-3})$. 
Assume that $\rho$ is $2$-thick. By Proposition 
\ref{prop:imp}, $\rho$ is irreducible. 
Suppose that $\wedge^2 \rho$ has a non-trivial $G$-invariant 
subspace $W \subset \wedge^2 V$ such that $2 \le \dim W \le 4$. 
Put $W_1 := W$ and $W_2 := W^{\perp} \subset \wedge^2 V$. 
Then $2 \le \dim W_2 \le 4$. 
By Lemma \ref{lemma:grass3}, $W_1$ and $W_2$ are realizable. 
Hence it follows from 
Proposition \ref{prop:condofm-thick} that $\rho$ is not $2$-thick.   
This is a contradiction. Therefore $\rho$ satisfies $(\ref{cond4-3})$ and 
we have $(\ref{cond4-2}) \Rightarrow (\ref{cond4-3})$. 

Finally, we show that $(\ref{cond4-3}) \Rightarrow (\ref{cond4-2})$. 
Assume that $\rho$ satisfies $(\ref{cond4-3})$. 
Suppose that $\rho$ is not $2$-thick.  
It follows from Proposition \ref{prop:condofm-thick} that    
there exist realizable invariant subspaces $W_1, W_2 \subset \wedge^2 V$ such 
that $W_1^{\perp} = W_2$. 
By Corollary \ref{cor:rge2} we have $\dim W_1 \ge 2$ and 
$\dim W_2 \ge 2$. 
Hence $W_1$ is a realizable invariant subspace such that $2 \le \dim W_1 \le 4$. 
This is a contradiction. Therefore $(\ref{cond4-3}) \Rightarrow (\ref{cond4-2})$. 
We have completed the proof. 
\hfill $\Box$

\bigskip 

Next, we deal with $5$-dimensional representations. 

%\begin{corollary}\label{cor:conj5} 
%Let $\rho :  G \to {\rm GL}(V)$ be an irreducible $5$-dimensional dimensional representation of a group $G$. 
%If the induced representation $\wedge^3 V$ has a 
%$G$-invariant realizable subspace $W$, then 
%$\dim W \ge {5 \choose 3} - 3(5-3) = 4$.  
%\end{corollary}
%
%{\it Proof}.
%Let us denote by $\chi : G \to {\rm GL}_1(k)$ the character induced by 
%$\wedge^5 \rho : G \to {\rm GL}(\wedge^5 V) \cong {\rm GL}_1(k)$.
%Let $V^{\ast} := \wedge^4 V$ be the dual $G$-module of $V$.  
%Then we can easily see that $\Phi : \wedge^3 V \otimes k_{\chi} \cong \wedge^2 V^{\ast}$ 
%as $G$-modules, where $k_{\chi}$ is the $1$-dimensional $G$-module defined by $\chi$. 
%Furthermore, this isomorphism $\Phi$ maps vectors associated to $3$-dimensional subspaces of $V$ 
%to vectors associated to $2$-dimensional subspaces of $V^{\ast}$. 
%If $\wedge^3 V$ has a non-trivial $G$-invariant realizable subspace $W$ with $0 < \dim W \le 3$, then 
%$\wedge^2 V^{\ast}$ also has a $G$-invariant realizable subspace $\Phi(W)$ 
%with $0 < \dim \Phi(W) \le 3$.  
%However, applying Lemma \ref{lemma:conj5} to the $G$-module $V^{\ast}$, 
%we see that $\dim \Phi(W) \ge 4$.  This is a contradiction. 
%Thus we have completed the proof.   
%\hfill $\Box$ 

\begin{proposition}\label{prop:condthick5} 
Let $V$ be a $5$-dimensional vector space over 
an algebraically closed field $k$.  
For a representation $\rho : G \to {\rm GL}(V)$,
the following are equivalent: 
\begin{enumerate}
\item\label{cond5-1} $\rho$ is thick.  
\item\label{cond5-2} $\rho$ is $2$-thick.  
\item\label{cond5-3} $\rho$ is irreducible and the induced representation 
$\wedge^2 \rho : G \to {\rm GL}(\wedge^2 V)$ has no non-trivial $G$-invariant 
subspace $W \subset \wedge^2 V$ with $4 \le \dim W \le 6$. 
\end{enumerate} 
\end{proposition}

{\it Proof}.
It is easy to check that $(\ref{cond5-1}) \Leftrightarrow (\ref{cond5-2})$ 
by Propositions~\ref{prop:dual} and \ref{prop:imp}. 
Let us show that $(\ref{cond5-1}) \Rightarrow (\ref{cond5-3})$. 
Assume that $\rho$ is thick.  
Then $\rho$ is irreducible by Proposition \ref{prop:imp}. 
Suppose that there exists a non-trivial $G$-invariant subspace $W \subset \wedge^2 V$ 
with $4 \le \dim W \le 6$. Put $W_1 := W$ and $W_2 := W^{\perp} \subset \wedge^3 V$.
By Lemma \ref{lemma:grass3}, $W_1$ and $W_2$ are realizable. 
Hence Proposition \ref{prop:condofm-thick} implies that $\rho$ is not thick. 
This is a contradiction. Therefore we see that $\rho$ satisfies $(\ref{cond5-3})$ and 
that $(\ref{cond5-1}) \Rightarrow (\ref{cond5-3})$. 

Finally, we show that $(\ref{cond5-3}) \Rightarrow (\ref{cond5-2})$. 
Assume that $\rho$ satisfies $(\ref{cond5-3})$. 
Suppose that $\rho$ is not $2$-thick. Then by Proposition \ref{prop:condofm-thick} 
there exist realizable subspaces 
$W_1 \subset \wedge^2 V$ and $W_2 \subset \wedge^3 V$ such that 
$W_1^{\perp} = W_2$. 
Proposition \ref{prop:conj5} says that $\dim W_1 \ge 4$ and $\dim W_2 \ge 4$. 
Hence $4 \le \dim W_1 \le 6$. This contradicts the assumption. 
Therefore $\rho$ is $2$-thick. 
\hfill $\Box$

\bigskip 

By using Proposition~\ref{prop:condthick5}, 
we have the following proposition.  

\begin{proposition}\label{prop:r-number2-5} 
$r(\wedge^2(5))=r(\wedge^3(5))=4$. 
\end{proposition} 

{\it Proof}. 
For a partition $\lambda = (\lambda_1, \ldots, \lambda_n)$ 
of $d$, we denote by $V_{\lambda}$ 
the irreducible representation of the symmetric group 
$S_d$ over ${\Bbb C}$ corresponding to $\lambda$. 
Let us consider the $5$-dimensional irreducible 
representations $V_{(3, 2)}$ and $V_{(2, 2, 1)}$ of 
$S_5$. By calculating characters, we see that 
$\wedge^2 V_{(3, 2)} = V_{(3, 1, 1)} \oplus V_{(2, 1, 1, 1)}$ and 
$\wedge^2 V_{(2, 2, 1)} = V_{(3, 1, 1)} \oplus V_{(2, 1, 1, 1)}$. 
Since $\dim V_{(3, 1, 1)} = 6$ and $\dim V_{(2, 1, 1, 1)} = 4$, 
$V_{(3, 2)}$ and $V_{(2, 2, 1)}$ are not $2$-thick by 
Proposition~\ref{prop:condthick5}. 
In particular, $\wedge^2 V_{(3, 2)}$ and 
$\wedge^2 V_{(2, 2, 1)}$ have $4$-dimensional $S_5$-invariant 
realizable subspaces $V_{(2, 1, 1, 1)}$, respectively. 
This implies that $r(\wedge^2(5)) = r(\wedge^3(5)) \le 4$. 
Using Proposition~\ref{prop:conj5}, we have 
$r(\wedge^2(5))=r(\wedge^3(5))=4$. 
\hfill $\Box$

\begin{remark}\label{remark:table}\rm 
By the calculations above, we have the following table on $r(\wedge^m(n))$ in the case 
$n \le 10$: 

\begin{center} 
\begin{tabular}{c|ccccccccccc}
$n \backslash m$ & $0$ & $1$ & $2$ & $3$ & $4$ & $5$ & $6$ & $7$ & $8$ & $9$ & $10$ \\  \hline  
$1$ & $1$ & $1$ & & & & & & & & & \\ 
$2$ & $1$ & $2$ & $1$ & & & & & & & & \\ 
$3$ & $1$ & $3$ & $3$ & $1$ & & & & & & & \\ 
$4$ & $1$ & $4$ & $2$ & $4$ & $1$ & & & & & & \\
$5$ & $1$ & $5$ & $4$ & $4$ & $5$ & $1$ & & & & & \\
$6$ & $1$ & $6$ & $3$ & $2$ & $3$ & $6$ & $1$ & & & & \\
$7$ & $1$ & $7$ & $A$ & $B$ & $B$ & $A$ & $7$ & $1$  & & & \\
$8$ & $1$ & $8$ & $4$ & $C$ & $2$ & $C$ & $4$ & $8$  & $1$ & & \\
$9$ & $1$ & $9$ & $D$ & $3$ & $E$ & $E$ & $3$ & $D$  & $9$ & $1$ & \\
$10$ & $1$ & $10$ & $5$ & $F$ & $G$ & $2$ & $G$ & $F$  & $5$ & $10$ & $1$ \\
\end{tabular} 
\end{center} 
Here we have not yet determined $A, B, C, D, E, F, G$. Until now, we can only know that   
$4 \le A \le7$, $3 \le B \le 7$, $3 \le C \le 8$, $5 \le D \le 9$, $3 \le E \le 9$, $4 \le F \le 10$, and $3 \le G \le 10$. 
\end{remark}

\bigskip 
 
 For $n \ge 6$, it is difficult to check 
whether a given $n$-dimensional representation 
is thick or not. 
In the rest of this section, we show 
some results on thickness and denseness.

\begin{lemma}\label{lemma:grphomthick} 
Let $\phi : G \to G'$ be a group homomorphism and 
$\rho : G' \to {\rm GL}(V)$ a finite-dimensional representation of $G'$. 
If $\rho$ is not $m$-thick, then neither is $\rho \circ \phi : G 
\to {\rm GL}(V)$.  
\end{lemma}

{\it Proof}. 
Suppose that $\rho \circ \phi$ is $m$-thick. 
Let $V_1$ and $V_2$ be subspaces of $V$ such that 
$\dim V_1 + \dim V_2 = \dim V$.  
Then there exists $g \in G$ such that 
$(\rho \circ \phi)(g) V_1 \oplus V_2 = V$. 
Putting $g' := \phi(g) \in G'$, we have $\rho(g')V_1 \oplus V_2 = V$. 
This implies $m$-thickness of $\rho$, which is a contradiction.  
Hence $\rho \circ \phi$ is not $m$-thick. 
\hfill $\Box$

\begin{proposition}\label{prop:wedge2kn}
Let $k$ be a field. Let $V := \wedge^2 k^n$ be the exterior of 
the standard representation $k^n$ of ${\rm GL}_n(k)$ with $n \ge 4$. 
Then $V$ is not $(n-1)$-thick as a representation of ${\rm GL}_n(k)$. 
\end{proposition}

{\it Proof}. 
Let $e_1, e_2, \ldots, e_n$ be the canonical basis of $k^n$. 
Let $W := \langle e_1 \wedge e_2, e_1 \wedge e_3, \ldots, e_1 \wedge e_n \rangle 
\subset \wedge^2 k^n$.  Note that the $(n-1)$-dimensional subspace $W$ is expressed as $e_1 \wedge k^n 
:= \{ e_1 \wedge v \mid v \in k^n \}$. 
For each $g \in {\rm GL}_n(k)$, $g W = \langle ge_1 \wedge ge_2, \ldots, 
ge_1 \wedge ge_n \rangle = (ge_1) \wedge k^n$. Put $ge_1 = a_1 e_1 + a_2 e_2 + \cdots + a_n e_n$. 
We see that $ge_1 \wedge e_1 \in g W$ and that   
$g e_1 \wedge e_1 = -a_2 e_1 \wedge e_2 - \cdots - a_n e_1 \wedge e_n \in W$.  
If $g e_1 \neq a_1 e_1$, then $0 \neq ge_1 \wedge e_1 \in W \cap gW$. 
If $g e_1 = a_1 e_1$, then $0 \neq ge_1 \wedge ge_2 = a_1 e_1 \wedge ge_2 \in W \cap g W$. 
Hence $W \cap gW \neq 0$. 
If we choose a subspace $W'$ of dimension $n(n-1)/2 - (n - 1)$ 
such that $W' \supset W$, 
then $gW \cap W' \neq 0$ for each $g \in W$.  
Hence $V$ is not $(n-1)$-thick. 
\hfill $\Box$

\begin{corollary}\label{cor:wedge2notthick}
Let $n \ge 4$.
For any $n$-dimensional representation $V$ of an arbitrary group $G$, 
the exterior representation $\wedge^2 V$ of $G$ is not $(n-1)$-thick. 
\end{corollary}

{\it Proof}. 
The statement follows from Lemma~\ref{lemma:grphomthick} and Proposition~\ref{prop:wedge2kn}.
\hfill $\Box$

\bigskip 

In the same way as Proposition~\ref{prop:wedge2kn}, we can prove the following proposition. This result was suggested by the referee. 

\begin{proposition}\label{prop:referee} 
Let $\rho  :  G \to {\rm GL}(V)$ be an $n$-dimensional representation of a group $G$ over a field $k$. Let $m \le \frac{1}{2} n$. Suppose that there exists an $m$-dimensional subspace $W$ of $V$ such that $(\rho(g) W) \cap W \neq 0$ for all $g \in G$. Then $\rho$ is not $m$-thick. 
\end{proposition} 

{\it Proof}. Since $n - m \ge 2m -m =m$, we can choose an $(n-m)$-dimensional 
subspace $W'$ of $V$ such that $W' \supseteq W$. For all $g \in G$, $(\rho(g)W) \cap W' \supseteq (\rho(g)W) \cap W' \cap W = (\rho(g)W) \cap W \neq 0$. 
Hence $\rho$ is not $m$-thick. 
\hfill $\Box$

\begin{remark}\rm 
Denseness and thickness are independent from absolute irreducibility. 
For example, the representation 
\[
\begin{array}{cccc} 
\rho : & {\mathbb R} & \to &  {\rm GL}(2, {\mathbb R}) \\
 & \theta & \mapsto & 
\left( 
\begin{array}{cc} 
\cos \theta & -\sin \theta \\
\sin \theta & \cos \theta \\
\end{array} 
\right)
\end{array} 
\]
is dense and thick, but not absolutely irreducible. 
Conversely, the representation $V = \wedge^2 {\mathbb C}^4$ of 
${\rm GL}(4, {\mathbb C})$ is not thick (and hence not dense) 
but absolutely irreducible. 
\end{remark}

The following proposition shows that there are many examples 
of representations which are not dense.

\begin{proposition} 
Let $n \ge 4$. 
Let $V$ be an $n$-dimensional irreducible representation of a group $G$. 
Assume that all irreducible representations of $G$ have dimension $\le n$. 
Then the representation $V$ of $G$ is not dense.   
\end{proposition}

{\it Proof}. 
The dimension of $\wedge^2 V$ is ${n \choose 2} (> n)$. Hence 
$\wedge^2 V$ can not be irreducible. This implies that $V$ is not dense.  
\hfill $\Box$

\begin{corollary}
Let $G$ be a finite group. 
Assume that $G$ has an irreducible representation of dimension $n$
with $n \ge 4$. 
Then $G$ has an irreducible representation which is not dense. 
\end{corollary} 

{\it Proof}.  
Let $n$ be the maximum of the dimensions of irreducible 
representations of $G$. Since $G$ is finite, there exists the maximum $n$. 
The assumption implies that $n \ge 4$. 
Let $V$ be an irreducible representation of $G$ of dimension $n$. 
The previous proposition shows that $V$ is not dense. 
\hfill $\Box$ 

\begin{remark}\rm 
Let $G$ be a group. Let $V$ be a representation of $G$ of dimension $n$ 
with $n \ge 4$ over ${\Bbb C}$. 
Then $\wedge^m V$ is not thick for $2 \le m \le n-2$. 
This fact will be proven in \cite{Classification}. 
\end{remark}

%%%%%%%%%%%%%%%%%%%%%%%%%%%%%%%%%%%%%%%%%%%%%%%%%%%%%%%%%%%%%%%%%%%%%

\section{Examples} 

In this section, we show several examples of representations for Lie groups. 
In \cite{Classification}, we will classify thick representations and dense 
representations for complex simple Lie groups. 
Here we show another approach to verify thickness and denseness. 

\subsection{Case: $G = {\rm GL}_{2}({\mathbb C})$} 
Let $a$ and $b$ be integers with $a \ge 0$. 
Let $V_{(a+b, b)}$ be the irreducible representation of $G = {\rm GL}_{2}({\mathbb C})$ 
with highest weight $(a+b, b)$. 
Set $\det^b := V_{(b, b)}$ for $b \in {\mathbb Z}$.  
Note that $V_{(a+b, b)} = \det^b \otimes V_{(a, 0)}$ and that 
$\dim V_{(a+b, b)} = a+1$. 

\begin{lemma}\label{lemma:gl2decomp} 
As representations of  ${\rm GL}_2({\Bbb C})$, we have 
\begin{eqnarray*} 
\wedge^2 V_{(a+b, b)} & = & \det {}^{2b} \otimes \left(
\sum_{k=1}^{[ \frac{a+1}{2} ]} \det {}^{2k-1} \otimes V_{(2a-4k+2, 0)}
\right) \\
 & = & \sum_{k=1}^{[ \frac{a+1}{2} ]} V_{(2a+2b-2k+1, 2b+2k-1)}. 
\end{eqnarray*} 
\end{lemma} 

{\it Proof}. Comparing characters, we can verify the statement. 
\hfill $\Box$

\begin{corollary}
For $a \ge 3$, the representation $V_{(a+b, b)}$ is not 2-dense.
In particular, it is not dense. 
\end{corollary}

{\it Proof}. Since $[ \frac{a+1}{2} ] \ge 2$ if $a \ge 3$, 
the representation $\wedge^2 V_{(a+b, b)}$ is not irreducible by 
Lemma \ref{lemma:gl2decomp}. 
Hence $V_{(a+b, b)}$ is not $2$-dense. 
\hfill $\Box$

\begin{corollary}
If $a = 3$ or $4$, then the representation $V_{(a+b, b)}$ is thick. 
\end{corollary}

{\it Proof}. 
When $a=3$, $\wedge^2 V_{(b+3, b)} = V_{(2b+5, 2b+1)} \oplus V_{(2b+3, 2b+3)}$ 
by Lemma \ref{lemma:gl2decomp}. Hence 
$\wedge^2 V_{(3+b, b)}$ has no ${\rm GL}_2({\Bbb C})$-invariant 
subspace $W$ such that $2 \le \dim W \le 4$ because 
$\dim V_{(2b+5, 2b+1)}=5$ and $\dim V_{(2b+3, 2b+3)} =1$.  
By Proposition \ref{prop:eq4}, the $4$-dimensional representation $V_{(b+3, b)}$ 
is thick.  When $a=4$, $\wedge^2 V_{(b+4, b)} = V_{(2b+7, 2b+1)} \oplus V_{(2b+5, 2b+3)}$ 
by Lemma \ref{lemma:gl2decomp}. Hence $\wedge^2 V_{(b+4, b)}$ 
has no ${\rm GL}_2({\Bbb C})$-invariant 
subspace $W$ such that $4 \le \dim W \le 6$ because 
$\dim V_{(2b+7, 2b+1)}=7$ and $\dim V_{(2b+5, 2b+3)} = 3$.  
By Proposition \ref{prop:condthick5}, the $5$-dimensional representation $V_{(b+4, b)}$ 
is thick.  
\hfill $\Box$ 

\begin{remark}\rm
When $a=1$ or $2$, the representation $V_{(a+b, b)}$ is dense. 
When $a \ge 3$, we can verify that $V_{(a+b, b)}$ is not dense, but thick  
by the classification of thick representations of simple Lie groups. 
Indeed, we will see that 
$S^m {\rm SL}_2$ is thick and not dense if $m \ge 3$ 
in \cite{Classification}. 
\end{remark}

\subsection{Case: $G = {\rm GL}_{n}({\mathbb C})$}

\begin{prop}
Let $V={\mathbb C}^n$ be the standard representation of ${\rm GL}_{n}({\mathbb C})$.  Then $V$ is dense. 
\end{prop}

{\it Proof}. 
This assertion follows from the irreducibility of 
$\wedge^i V $ for $1 \le i \le n-1 $. 
\hfill $\Box$

\bigskip 

For the standard representation $V = {\Bbb C}^n$ of 
${\rm GL}_n({\Bbb C})$, let us discuss thickness and 
denseness of ${\rm S}^2 V$ and $\wedge^2 V$.

\begin{lemma}\label{lemma:numberofpartitions} 
Put  $P_n(x)=\prod_{i=1}^n (1+x^i) = 
(1+x)(1+x^2)\cdots (1+x^n)$.
Let $a_{i}$ be the coefficient of $x^{i} $ in $ P_n(x)$.  
Then  if $n \ge 3$, then $a_{i} \ge 1$ for any $0 \le i \le \frac{n(n+1)}{2} $ and $a_{i} \ge 2$ for any $3 \le i \le \frac{n(n+1)}{2}-3$.
\end{lemma}
{\it Proof}. 
Let us prove the statement by induction on $n$. 
When $n=3$, $P_3(x)=1+x+x^2+2x^3+x^4+x^5+x^6$ and hence the statement holds. 
Suppose that the statement is true for $n$. 
Since 
\begin{multline*} 
P_{n+1}(x)  =   P_n(x)(1+x^{n+1}) \\ 
  =   a_0 +a_1x + a_2 x^2 + \cdots  + a_n x^n + (a_{n+1}+a_0)x^{n+1} 
+(a_{n+2}+a_1)x^{n+2}  \\ 
+ \cdots   + (a_{n(n+1)/2} + a_{(n+1)(n-2)/2})x^{n(n+1)/2} +a_{n(n-1)/2}x^{(n^2+n+2)/2} \\ 
+ \cdots + a_{n(n+1)/2} x^{(n+1)(n+2)/2},  
\end{multline*}  
the statement is also true for $n+1$. This completes the proof. 
\hfill $\Box$

\begin{prop}\label{prop:6.7}
Let $V={\mathbb C}^n$ be the standard representation of ${\rm GL}_{n}({\mathbb C})$ $(n\ge3)$. 
Then the second symmetric tensor ${\rm S}^2 V$ is irreducible, but not $m$-$thick$ for $3 \le m
 \le \frac{n(n+1)}{2} - 3$. 
\end{prop}

{\it Proof}. 
It is well-known that ${\rm S}^2 V$ is irreducible. 
By \cite[Theorem~4.4.2]{Howe}, the number of irreducible 
components of  $ \wedge^m ({\rm S}^2 V)$ is equal to the 
the number of partitions of $m$ into distinct parts of 
size at most $n$. 
This number is equal to the coefficient $a_m$ of $x^m$ in $P_{n}(x)$ 
in Lemma~\ref{lemma:numberofpartitions}. 
Since $a_m \ge 2$ by Lemma~\ref{lemma:numberofpartitions}, $\wedge^m ({\rm S}^2 V)$ is not irreducible  for any $3 \le m
\le \frac{n(n+1)}{2} - 3$. 
We see that irreducible components of $\wedge^m ({\rm S}^2 V)$ are all 
realizable by the proof of \cite[Theorem~4.4.2]{Howe}. Hence 
Proposition \ref{prop:condofm-thick} implies that $\wedge^m ({\rm S}^2 V)$ is not $m$-thick for any $3 \le m
\le \frac{n(n+1)}{2} - 3$. 
\hfill $\Box$

\bigskip

Proposition \ref{prop:wedge2kn} shows that $\wedge^2 {\Bbb C}^n$ 
is not $(n-1)$-thick for the standard representation ${\Bbb C}^n$ of 
${\rm GL}_n({\Bbb C})$ for $n \ge 4$.  Moreover, we have the following proposition.

\begin{prop}\label{prop:6.8}
Let $V={\mathbb C}^n$ be the standard representation of ${\rm GL}_{n}({\mathbb C})$ $(n\ge4)$. 
Then the second alternating tensor $\wedge^2 V$ is irreducible, but not $m$-$thick$ for $3 \le m
 \le \frac{n(n-1)}{2} - 3$. 
\end{prop}
{\it Proof}. 
It is well-known that $\wedge^2 V$ is irreducible. 
By \cite[Theorem~4.4.4]{Howe}, the number of irreducible 
components of  $ \wedge^m (\wedge^2 V)$ is equal to the 
the number of partitions of $m$ into distinct parts of 
size at most $n-1$. 
This number is equal to the coefficient $a_m$ of $x^m$ in $P_{n-1}(x)$ 
in Lemma~\ref{lemma:numberofpartitions}. 
Since $a_m \ge 2$ by Lemma~\ref{lemma:numberofpartitions}, $\wedge^m (\wedge^2 V)$ is not irreducible  for any $3 \le m
\le \frac{n(n-1)}{2} - 3$. 
We see that irreducible components of $\wedge^m (\wedge^2 V)$ are all 
realizable by the proof of \cite[Theorem~4.4.4]{Howe}. Hence 
Proposition \ref{prop:condofm-thick} implies that $\wedge^m (\wedge^2 V)$ is not $m$-thick for any $3 \le m
\le \frac{n(n-1)}{2} - 3$. 
\hfill $\Box$

\bigskip 

Then from Propositions \ref{prop:6.7} and \ref{prop:6.8},  we have the following corollary.

\begin{corollary} 
Let $V={\mathbb C}^n$ be an $n$-dimensional representation of any group $G$. 
If $n \ge 3$, then the second symmetric tensor ${\rm S}^2 V$ 
is not thick. If $n \ge 4$, then the second alternating tensor 
$\wedge^2 V$ is not $thick$. 
\end{corollary}
 
{\it Proof}. 
Using Lemma \ref{lemma:grphomthick}, we can prove the statement  
(the latter part has been proved in Corollary \ref{cor:wedge2notthick}). 
\hfill $\Box$

\subsection{  Case: $G = {\rm SO}_{n}({\mathbb C})$}
\begin{prop}\label{prop:so2ndensethick} 
Let $V$ be the standard representation of $G = {\rm SO}_{2n}({\mathbb C})$.  Then $V$ is $m$-dense for each $0 < m < 2n$ with $m\neq n$, but not $n$-thick.
\end{prop}

{\it Proof}. 
The first assertion follows from the irreducibility of $\wedge^i V$ for $1 \le i \le n-1$. 
The proof of \cite[Theorem 19.2]{FH} shows that the $n$-th alternating tensor $\wedge^n V$ has exactly two irreducible factors and they are realizable. Then by Proposition \ref{prop:condofm-thick} the representation $V$ is not $n$-thick. 
\hfill $\Box$

\begin{prop}
Let $V$ be the standard representation of $G = {\rm SO}_{2n+1}( {\mathbb C})$. Then $V$ is dense.  
\end{prop}

{\it Proof}. 
The $m$-th alternating tensor $\wedge^m V$ is irreducible for each $0 < m < 2n+1$ 
(for example, see \cite[Theorem~19.14]{FH}).  This implies the statement. 
\hfill $\Box$

\subsection{  Case: $G = {\rm Sp}_{2n}({\mathbb C})$}

Let $V$ be a $2n$-dimensional complex vector space, $\{ e_1,e_2,\allowbreak \ldots,e_{2n} \}$ a basis for $V$, and $\{ e^{*}_1,e^{*}_2,\ldots,e^{*}_{2n} \}$ its dual basis for the dual vector space $V^{*}$.
We use a non-degenerate skew-symmetric bilinear form $\omega =\displaystyle{\sum ^{n}_{i=1}}e^{*}_{i}\wedge e^{*}_{n+i}$ and the corresponding symplectic Lie group 
${\rm Sp}_{2n}({\mathbb C})$. 
Then we have a contraction map by $\omega$:
 $$f_{m}:\wedge^{m}V\rightarrow \wedge^{m-2}V.$$
 If $m \le n $, $\mathrm{Ker} f_{m}$ is the $m$-th fundamental representation of 
${\rm Sp}_{2n}({\mathbb C})$. In particular, $\mathrm{Ker} f_{m}$ is irreducible. 
We have the isotropic Grassmann variety of isotropic subspaces of dimension $m$ as a unique minimal closed orbit in the projective space $\bold{P}(\mathrm{Ker}f_{m})$. 
Since $\mathrm{Ker} f_{m}$ contains $\wedge^{m} L$ for any isotropic $m$-dimension subspace $L \subset V$, $\mathrm{Ker} f_{m}$ is realizable. 
For details see  \cite{FH}.

The following lemma is well-known. 
\begin{lemma}\label{lemma:basis}
Let $(V,\omega)$ be a $2n$-dimensional symplectic vector space and $W\subset V$ a subspace. Then there is a basis $\{ v_1,v_2,\ldots , v_{2n} \}$ of $V$ such that 
$\omega (v_i,v_{n+i})=1$, $\omega (v_i,v_j)=0\  \mathrm{if} \  j\neq i\pm n$, and for some non-negative integers $l,k$
$$W=\langle v_1,\ldots ,v_{k},v_{k+1},\ldots ,v_{n-l},v_{n+1},\ldots , v_{n+k}\rangle.$$ 
\end{lemma}

\begin{lemma}\label{lemma:lagrangian}
Let $(V,\omega)$ be a $2n$-dimensional symplectic vector space and $W\subset V$ a subspace of codimension $i$ $(i\le n)$. Then there is a Lagrangian subspace $L\subset V$ such that $L + W= V$. 
\end{lemma} 
{\it Proof}. 
It is enough to prove the case of $i=n$.
By Lemma \ref{lemma:basis}, 
there is a symplectic basis $\{ v_1,v_2,\ldots , v_{2n} \}$ of $V$ such that 
for some non-negative integer $k\le \frac{n}{2}$
$$W=\langle v_1,\ldots ,v_{k},v_{k+1},\ldots ,v_{n-k},v_{n+1},\ldots , v_{n+k}\rangle.$$
In the case of $k=0$, $W=\langle v_1,\ldots ,v_{n}\rangle$. Then $L=\langle  v_{n+1},\ldots , v_{2n}\rangle$ satisfies the condition $L + W= V$. 
In the case of $1 \le k\le \frac{n}{2}$, we put as following,
$$L=\langle v_{n+k+1},\ldots ,v_{n+k+i},\ldots ,v_{2n-k}, \ \ \ \ \ \ \ \ \ \ \ \ \ \ \ \ \ \ \ \ \ \ \ \ \ \ \ \ \ \ \ \ \ \ \ \ \ \ \ \ \ \ \ \ \ \ \ \ $$
$$ v_{n-k+1}+v_{n+1},\ldots ,v_{n-k+i}+v_{n+i},\ldots ,v_{n}+v_{n+k}, $$
$$\ \ \ \ \ \ \ \ \ \ \ \ \ \ \ \ \ \ \ \ \ \ \ \ \ v_{2n-k+1}+v_{1},\ldots ,v_{2n-k+i}+v_{i},\ldots ,v_{2n}+v_{k} \rangle.$$
Then $L$ is a Lagrangian subspace and satisfies the condition $L + W= V$. 
\hfill $\Box$

\begin{lemma}\label{lemma:isotropic}
Let $(V,\omega)$ be a $2n$-dimensional symplectic vector space and $W\subset V$ a subspace of codimension $i$ $(i\le n)$. Then there is an isotropic subspace $U\subset V$ of dimension $i$ such that $U \cap W= \{ 0\}$. 
\end{lemma} 
{\it Proof}. 
By Lemma \ref{lemma:lagrangian}, there is a Lagrangian subspace $L\subset V$ such that $L + W= V$. 
Since the dimension of $L\cap W$ is $n-i$, there is a subspace $U\subset L$ such that the dimension of $U$ is $i$ and  $U \cap W= \{ 0\}$.
Since $U$ is a subspace of a Lagrangian subspace $L$, $U$ is an isotropic subspace.  
\hfill $\Box$ 

\bigskip

Then we have the following proposition. 
\begin{prop}\label{prop:ker}
Let $(V,\omega)$ be the standard representation of %symplectic group 
${\rm Sp}_{2n}({\mathbb C})$.
For each $1 < m \le n$, $(\mathrm{Ker}f_{m})^{\perp} \subset \wedge^{2n-m}V$ is not realizable.
\end{prop}
{\it Proof}. 
If $(\mathrm{Ker}f_{m})^{\perp}$ is realizable, there is a subspace $W\subset V$  of codimension $m$ such that $\wedge^{2n-m}W \in (\mathrm{Ker}f_{m})^{\perp}$. Then by Lemma \ref{lemma:isotropic} we have an isotropic subspace $U\subset V$ of dimension $m$ such that $U \cap W= \{ 0\}$. Because $\mathrm{Ker}f_{m}$ contains $\wedge^{m} L$ for any isotropic subspace $L \subset V$ of dimension $m$, we have $\wedge^{m} U \in \mathrm{Ker}f_{m}$. But we have $(\wedge^{2n-m} W )\wedge (\wedge^{m}U)\neq 0$. This is a contradiction.   
\hfill $\Box$ 
\bigskip

By the ${\rm SL}_{2n}({\mathbb C})$-equivariant canonical pairing $\wedge^{2n-k}V \times \wedge^{k}V \rightarrow \wedge^{2n}V \cong \mathbb C$, we have the ${\rm SL}_{2n}({\mathbb C})$-equivariant isomorphism 
$$\wedge^{2n-k}V \rightarrow (\wedge^{k}V)^{*} \cong \wedge^{k}V^{*}. $$

Moreover  by the correspondence $e^{*}_{i_{1}}\wedge \cdots \wedge e^{*}_{i_{k}}  \mapsto e_{i_{1}}\wedge \cdots \wedge e_{i_{k}} $ we have the isomorphism $\wedge^{k}V^{*} \rightarrow \wedge^{k}V$ as vector spaces.
The difference between these vector spaces as ${\rm SL}_{2n}({\mathbb C})$-modules is  described by the outer automorphism 
\[ 
\begin{array}{ccccc}
\sigma & : & {\rm SL}_{2n}({\mathbb C})   & \rightarrow & 
{\rm SL}_{2n}({\mathbb C})  \\
       &   & g                          & \mapsto     & {}^{t}g^{-1} .
\end{array}
\]

Then we obtain the isomorphism $\phi$ as ${\rm SL}_{2n}({\mathbb C})$-modules up to the outer automorphism $\sigma$, that is,  
$$\phi : \wedge^{2n-k}V \rightarrow (\wedge^{k}V)^{*} \cong \wedge^{k}V^{*} \rightarrow \wedge^{k} V.$$
Thereby $\phi$ induces the isomorphism $\overline{\phi}$ as follows:  
 
\[ 
\begin{array}{ccccc}
\overline{\phi} & : & {\mathbb P}(\wedge^{2n-k} V) & \rightarrow & {\mathbb P}(\wedge^{k} V)  \\
                &   & \cup                          &             & \cup              \\
                &   & {\rm Gr}(2n-k,V)                 &  \cong     &  {\rm Gr}(k,V).                
\end{array}
\]

Since $\sigma ({\rm Sp}_{2n}({\mathbb C}))={\rm Sp}_{2n}({\mathbb C})$ 
and $\sigma$ is an inner automorphism of ${\rm Sp}_{2n}({\mathbb C})$,
$\phi$ gives an 
isomorphism between $\wedge^{2n-k}V$ and $\wedge^{k}V$ as 
${\rm Sp}_{2n}({\mathbb C})$-modules. 
When we consider $\wedge^{k}V$ as a ${\rm Sp}_{2n}({\mathbb C})$-module, it is 
well-known that each irreducible representation of ${\rm Sp}_{2n}({\mathbb C})$ occurs at most once in an irreducible decomposition   
of $\wedge^{k}V$ (see  \cite[Chap.17]{FH}). Then  we have several irreducible ${\rm Sp}_{2n}({\mathbb C})$-invariant subspaces $\{ W_{i} \}_{i=1,\ldots ,s}$ in $\wedge^{k}V$ such that we have a unique irreducible decomposition $\wedge^{k}V= W_{1}\oplus W_{2} \oplus  \cdots \oplus W_{s} $, and $W_{i} \cong W_{j}$ if and only if $i=j$.   
Since there exists some number $i$ such that $W_{i}= {\mathrm{Ker}f_{k}}$, from now we put $W_{1}= {\mathrm{Ker}f_{k}}$.
Therefore under the isomorphism $\phi$ we can obtain the unique irreducible decomposition of $\wedge^{2n-k}V$. Namely if we put $W'_{i}:={\phi}^{-1}(W_{i})$, $\{ W'_{i} \}_{i=1,\ldots ,s}$ are ${\rm Sp}_{2n}({\mathbb C})$-invariant subspaces in $\wedge^{2n-k}V$ such that we have the unique irreducible decomposition $\wedge^{2n-k}V= W'_{1}\oplus W'_{2} \oplus \cdots \oplus W'_{s} $, and $W'_{i} \cong W'_{j}$ if and only if $i=j$.   
By the above construction we have the following lemma.

\begin{lemma}\label{lemma:realSp}
For any subset $\{ j_{1},\ldots ,j_{l}  \}  \subset \{1,2,\ldots ,s \}$, the subset ${\mathbb P}(W_{j_1}\oplus \cdots \oplus W_{j_{l}}) \cap {\mathrm Gr}(k, V)$ is empty if and only if the subset ${\mathbb P}(W'_{j_1}\oplus \cdots \oplus W'_{j_{l}}) \cap {\mathrm Gr}(2n-k, V)$ is empty. 
%{\color{red} In particular, $W_{j_1}\oplus \cdots \oplus W_{j_{l}} \subset 
%\wedge^k V$ is realizable if and only if $W^{'}_{j_1}\oplus \cdots \oplus W^{'}_{j_{l}} 
%\subset  \wedge^{2n-k} V$ is realizable. }
\end{lemma} 

\begin{prop}\label{Corollary:realSp}
For any subset $\{ j_{1},\ldots ,j_{l}  \}  \subset \{1,2,\ldots ,s \}$, the following are equivalent:
\begin{enumerate}
\item\label{condition1} $W_{j_1}\oplus \cdots \oplus W_{j_{l}}$ is a realizable subspace of $\wedge^k V$. 

\item\label{condition2} $W'_{j_1}\oplus \cdots \oplus W'_{j_{l}}$ is a realizable subspace of $\wedge^{2n-k} V$.

\item\label{condition3} There is some $m \in \{ 1, \ldots ,l\}$ such that $j_{m}=1.$
\end{enumerate}

\end{prop} 
{\it Proof}. 
Lemma \ref{lemma:realSp} 
shows that (\ref{condition1}) and (\ref{condition2}) are equivalent. 
Note that $W^{\ast} \cong W$  
for any ${\rm Sp}_{2n}({\mathbb C})$-modules $W$.  
For the perfect paring $\wedge^k V \times \wedge^{2n-k} V \to 
\wedge^{2n} V \cong k$, we see that $({\mathrm{Ker}f_{k}})^{\perp} = W'_{2}\oplus W'_{3} \oplus \cdots \oplus W'_{s}$. 
Indeed, any ${\rm Sp}_{2n}({\Bbb C})$-homomorphism 
$\phi : W_1 = \mathrm{Ker}f_{k} \to (W'_{2}\oplus W'_{3} \oplus \cdots \oplus W'_{s})^{\ast} \cong 
W'_{2}\oplus W'_{3} \oplus \cdots \oplus W'_{s}$ is zero.  
By Corollary \ref{cor:perfperp}, we have $({\mathrm{Ker}f_{k}})^{\perp} 
\supseteq W'_{2}\oplus W'_{3} \oplus 
\cdots \oplus W'_{s}$. 
Since $W'_{1} \cong W_1$ is irreducible, 
 $({\mathrm{Ker}f_{k}})^{\perp} = W'_{2}\oplus W'_{3} \oplus \cdots \oplus W'_{s}$. 
Then Proposition \ref{prop:ker} shows  that  (\ref{condition2}) and  (\ref{condition3}) are equivalent. 
\hfill $\Box$ 

\bigskip 

Then we have the following proposition.
\begin{prop}
The standard representation of ${\rm Sp}_{2n}({\mathbb C})$ is thick, but not $m$-dense for each $1<m<2n-1$.
\end{prop}
{\it Proof}. 
Since each irreducible representation occurs at most once in $\wedge^{k}V$, 
for any invariant subspace $U\subset \wedge^{k}V$ 
there is a subset $\{ i_{1},\ldots ,i_{\alpha}  \}  \subset \{1,2,\ldots ,s \}$ 
such that $U=W_{i_1}\oplus \cdots \oplus W_{i_{\alpha}}$. Similarly for $U^{\perp}$ there is a subset $\{ j_{1},\ldots ,j_{\beta} \}  \subset \{1,2,\ldots ,s \}$ 
such that $U^{\perp}
=W'_{j_1}\oplus \cdots \oplus W'_{j_{\beta}}$.
Since $({\mathrm{Ker}f_{k}})^{\perp} = W'_{2}\oplus W'_{3} 
\oplus \cdots \oplus W'_{s}$, 
$1 \in \{ i_{1},\ldots ,i_{\alpha}  \}$ if and only if $1 \notin  
\{ j_{1},\ldots ,j_{\beta}  \}$. 
By Proposition \ref{Corollary:realSp}, it is impossible that 
both $U$ and $U^{\perp}$ are realizable.  
This implies that $V$ is thick. 
Since it is well-known that $\wedge^m V$ is not irreducible for each $1<m<2n-1$,  
$V$ is not $m$-dense for $1<m<2n-1$.  
\hfill $\Box$

%\section{Conjectures} 
%\begin{conjecture} 
%For an $n$-dimensional representation 
%$\rho : G \to {\rm GL}(V)$, 
%is it true that $(m+1)$-thick $($or $(m+1)$-dense$)$ implies $m$-thick 
%$($or $m$-dense, respectively$)$ for 
%$0 < m \le n/2?$ 
%Is it true that $(i, j)$-thick $($or $(i, j)$-dense$)$ implies $(i-1, j+1)$-thick 
%$($or $(i-1, j+1)$-dense, respectively$)$ for 
%$0 < i \le j?$
%\end{conjecture} 
%
%\begin{conjecture}
%$(i-1, j)$-thick and $(i, j-1)$-thick implies $(i, j)$-thick?
%\end{conjecture} 

\end{document}